\documentclass[10pt]{article}

\usepackage{amsfonts} 
\usepackage{amsmath,hhline,latexsym,mathrsfs}
\usepackage{amssymb}
\usepackage{relsize}
\DeclareMathAlphabet\mathbfcal{OMS}{cmsy}{b}{n}

\usepackage[latin1]{inputenc}

\usepackage{hyperref}

\usepackage{color}
\usepackage{graphicx}
\usepackage{comment}

\newtheorem{defi}{Definition}[section]
\newtheorem{theorem}{Theorem}
\newtheorem{prop}{Proposition}
\newtheorem{lemma}{Lemma}

\newtheorem{remark}{Remark}

\newcommand{\R}{\mathbb{R}}
\newcommand{\N}{\mathbb{N}}
\newcommand{\A}{\mathcal{A}}

\newcommand{\ph}{\varphi}

\begin{document}

\title{Approximation of exact controls for semi-linear 1D wave equations using a least-squares approach}

\author{
\textsc{Arnaud M\"unch}\thanks{Laboratoire de Math\'ematiques Blaise Pascal, Universit\'e Clermont Auvergne, 
UMR CNRS 6620, Campus universitaire des C\'ezeaux, 3, place Vasarely, 63178, Aubi\`ere, France. E-mail: {\tt arnaud.munch@uca.fr.}}
\and
\textsc{Emmanuel Tr\'elat}\thanks{Laboratoire Jacques-Louis Lions, Sorbonne Universit\'e, 4 place Jussieu, 75012 Paris, France. E-mail: {\tt emmanuel.trelat@sorbonne-universite.fr.}}
}

\maketitle

\begin{abstract}
The exact distributed controllability of the semilinear wave equation $y_{tt}-y_{xx} + g(y)=f \,1_{\omega}$, assuming that $g$ satisfies the growth condition $\vert g(s)\vert /(\vert s\vert \log^{2}(\vert s\vert))\rightarrow 0$ as $\vert s\vert \rightarrow \infty$ and that $g^\prime\in L^\infty_{loc}(\mathbb{R})$ has been obtained by Zuazua in the nineties. The proof based on a Leray-Schauder fixed point argument makes use  of precise estimates of the observability constant for a linearized wave equation. It does not provide however an explicit construction of a null control.   
Assuming that $g^\prime\in L^\infty_{loc}(\mathbb{R})$, that $\sup_{a,b\in \mathbb{R},a\neq b} \vert g^\prime(a)-g^{\prime}(b)\vert/\vert a-b\vert^r<\infty $ for some $r\in (0,1]$ and that $g^\prime$ satisfies the growth condition $\vert g^\prime(s)\vert/\log^{2}(\vert s\vert)\rightarrow 0$ as $\vert s\vert \rightarrow \infty$, we construct an explicit sequence converging strongly to a null control for the solution of the semilinear equation. The method, based on a least-squares approach guarantees the convergence whatever the initial element of the sequence may be. In particular, after a finite number of iterations, the convergence is super linear with rate $1+r$. This general method provides a constructive proof of the exact controllability for the semilinear wave equation. 
\end{abstract}



\section{Introduction}

 Let $\Omega:=(0,1)$ and $\omega$ be a non empty open set of $\Omega$.  Let $T>0$ and denote $Q_T:=\Omega\times (0,T)$, $q_T:=\omega\times (0,T)$ and $\Sigma_T:=\partial\Omega\times (0,T)$.
 Let $g:\mathbb{R}\to \mathbb{R}$ be a continuous function. We consider the semilinear wave equation
\begin{equation}
\label{eq:wave-NL}
\left\{
\begin{aligned}
& y_{tt} - y_{xx} +  g(y)= f 1_{\omega}, & \textrm{in}\,\, Q_T,\\
& y=0, & \textrm{on}\,\, \Sigma_T, \\
& (y(\cdot,0),y_t(\cdot,0))=(u_0,u_1), & \textrm{in}\,\, \Omega,
\end{aligned}
\right.
\end{equation}
where $(u_0,u_1)\in \boldsymbol{V}:=H^1_0(\Omega)\times L^2(\Omega)$ is the initial state of $y$ and $f\in L^2(q_T)$ is a {\it control} function. We assume that there exists a positive constant $C$ such that
\begin{equation}
\label{assumption_g}
\vert g(s)\vert \leq C (1+\vert s\vert)\log^2(1+\vert s\vert), \quad \forall s\in\mathbb{R},
\end{equation}
so that (\ref{eq:wave-NL}) has a unique global weak solution in $C([0,T]; H_0^1(\Omega)) \cap C^1([0,T]; L^2(\Omega))$ (see \cite{CazenaveHaraux1980}).

The exact controllability for \eqref{eq:wave-NL} is formulated as follows: for any $(u_0,u_1), (z_0,z_1)\in \boldsymbol{V}$, find a control function $f\in L^2(q_T)$ such that the weak solution of \eqref{eq:wave-NL} satisfies $(y(\cdot,T),y_t(\cdot,T))=(z_0,z_1)$ at time $T>0$. 
Assuming $T$ large enough and a growth condition on the nonlinearity $g$ at infinity, this problem has been solved in \cite{zuazua93}.
\begin{theorem}\label{zuazuaTH}\cite{zuazua93} Assume that $\omega=(l_1,l_2)$ with $0\leq l_1<l_2\leq 1$. Assume $T>2\max(l_1,1-l_2)$. 
There exists a $\beta>0$ such that if the function $g\in C^1(\mathbb{R})$ satisfies
\begin{equation}\label{asymptotic_behavior}
\limsup_{\vert s\vert \to \infty} \frac{\vert g(s)\vert }{\vert s\vert \log^2\vert s\vert} < \beta,
\end{equation}
then for all $(y_0,y_1), (z_0,z_1)\in \boldsymbol{V}$, there exists control 
$f\in L^2(q_T)$ such that the solution of \eqref{eq:wave-NL} satisfies $(y(\cdot,T),y_t(\cdot,T))=(z_0,z_1)$.
\end{theorem}


The constant $\beta$ depends on $T$ and $\Omega$  but is independent of the data $(u_0,u_1)$, $(z_0,z_1)$.
Moreover, it is proved in \cite{zuazua93} that if $g$ behaves at infinity like $-s\log^p(\vert s\vert)$ with $p>2$, then due to blow up phenomena, the system is not exactly controllable in any time $T>0$. Later on, Theorem \ref{zuazuaTH} has been improved in \cite{Cannarsa_Loreti_Komornik2002} relaxing the condition \eqref{asymptotic_behavior} by 
$$
\limsup_{\vert s\vert\to \infty} \biggl\vert \int_0^s g(r)dr\biggl\vert  \biggl[\vert s\vert \prod_{k=1}^\infty \log_k(e_k+s^2)\biggr]^{-2}<\infty
$$
where $\log_k$ denotes the iterate logarithm function and $e_k$ the real number such that $\log_k(e_k)=1$. This growth condition is sharp since solution of \eqref{eq:wave-NL} may blow up whenever $g$ grow faster at infinity and $g$ has the bad sign. 
The multidimensional case for which $\Omega$  is a bounded domain of $\mathbb{R}^d$, $d>1$ with a $C^{1,1}$ boundary has been addressed in \cite{Li_Zhang_2000}. Assuming that the support $\omega$ of the control function is a neighborhood of $\partial\Omega$ and that $T>\textrm{diam}(\Omega\backslash\omega)$, the exact controllability of (\ref{eq:wave-NL})  is proved assuming the growth condition 
$\limsup_{\vert s\vert \to \infty} \frac{\vert g(s)\vert }{\vert s\vert\sqrt{\log\vert s\vert}} < \infty$. For control domains $\omega$ satisfying the classical multiplier method of Lions \cite{JLL88}, the exact controllability has been proved in \cite{Zhang_2000} assuming $g$ globally Lipschitz continuous. We also mention \cite{coron-trelat-wave-06} where a positive boundary controllability result is proved for a specific class of initial and final data and $T$ large enough.

The proof given in \cite{zuazua93} is based on a fixed point argument introduced in \cite{zuazua91}  that reduces the exact controllability problem to the obtention of suitable \textit{a priori} estimates for the linearized wave equation with a potential. More precisely, it is shown that the operator $K:L^\infty(Q_T)\to L^\infty(Q_T)$, where $y_\xi:=K(\xi)$ is a controlled solution through the control function $f_{\xi}$ of the linear boundary value problem 
\begin{equation}
\label{NL_z}
\left\{
\begin{aligned}
& y_{\xi,tt} - y_{\xi,xx} +  y_\xi \,\widehat{g}(\xi)= -g(0)+f_\xi 1_{\omega}, &\textrm{in}\,\, Q_T,\\
& y_\xi=0, &\textrm{on}\,\, \Sigma_T, \\
& (y_\xi(\cdot,0),y_{\xi,t}(\cdot,0))=(u_0,u_1), &\textrm{in}\,\, \Omega,
\end{aligned}
\right.
\qquad
\widehat{g}(s):=
\left\{ 
\begin{aligned}
& \frac{g(s)-g(0)}{s} & s\neq 0,\\
& g^{\prime}(0) & s=0
\end{aligned}
\right., 
\end{equation}
satisfying $(y_{\xi}(\cdot,T),y_{\xi,t}(\cdot,T))=(z_0,z_1)$ possesses a fixed point. The control $f_{\xi}$ is chosen in \cite{zuazua93} as the one of minimal $L^2(q_T)$-norm. The existence of a fixed point for the compact operator $K$ is obtained by using the Leray-Schauder's degree theorem. Precisely, it is shown that if $\beta$ is small enough, then there exists a constant $M=M(\Vert u_0,u_1\Vert_{\boldsymbol{V}}, \Vert z_0,z_1\Vert_{\boldsymbol{V}}$) such that $K$ maps the ball $B_\infty(0,M)$ into itself.

 The main goal of this work is to determine an approximation of the controllability problem associated to~(\ref{eq:wave-NL}), that is to construct an explicit sequence $(f_k)_{k\in \mathbb{N}}$ converging strongly toward an exact control for~(\ref{eq:wave-NL}). A natural strategy is to take advantage of the method used in \cite{zuazua93} and consider the Picard iterates $(y_k)_{k>0}$ associated with the operator $K$ defined by $y_{k+1}=K(y_k)$, $k\geq 0$ initialized with any element $y_0\in L^{\infty}(Q_T)$. The sequence of controls is then $(f_k)_{k\in \mathbb{N}}$ so that $f_{k+1}\in L^2(q_T)$ is the control of minimal $L^2(q_T)$-norm for $y_{k+1}$ solution of 
\begin{equation}
\label{NL_z_k}
\left\{
\begin{aligned}
& y_{k+1,tt} - y_{k+1,xx} +  y_{k+1} \,\widehat{g}(y_k)= -g(0)+f_{k+1} 1_{\omega}, & \textrm{in}\quad Q_T,\\
& y_{k+1}=0, & \textrm{on}\,\, \Sigma_T, \\
& (y_{k+1}(\cdot,0),y_{k+1,t}(\cdot,0))=(y_0,y_1), & \textrm{in}\,\, \Omega.
\end{aligned}
\right. 
\end{equation}
Such strategy usually fails since the operator $K$ is a priori not strictly contracting including for globally Lipschitz continuous function $g$. We refer to \cite{EFC-AM-2012} exhibiting lack of convergence in parabolic cases when such method is used. 
%
%
As is also usual for nonlinear problems, we may employ a Newton type method in order to find a zero of the mapping $\widetilde{F}: Y \mapsto W$ defined by 
   \begin{equation}\label{def-F}
\widetilde{F}(y,f):= \biggl(y_{tt} - y_{xx}  + g(y) - f 1_{\omega}, y(\cdot\,,0) - u_0, y_t(\cdot\,,0) - u_1, y(\cdot\,,T) -z_0, y_t(\cdot\,,T) - z_1\biggr) 
   \end{equation}
   for some appropriates Hilbert spaces $Y$ and $W$ (see below). 
%
%
Assuming $g\in C^1(\mathbb{R})$ so that $\widetilde{F}\in C^1(Y,W)$, the Newton iterative method for $\widetilde{F}$ reads as follows: given $(y_0,f_0)$ in $Y$, define the sequence $(y_k,f_k)_{k\in \mathbb{N}}$ iteratively as follows $(y_{k+1},f_{k+1})=(y_k,f_k)-(Y_k,F_k)$ where $F_k$ is a control for $Y_k$ solution of  
\begin{equation}
\label{Newton-nn}
\left\{
\begin{aligned}
&Y_{k,tt}-Y_{k,xx} +  g^{\prime}(y_k)\,Y_{k} = F_{k}\, 1_{\omega} + y_{k,tt} - y_{k,xx}  + g(y_k) - f_k 1_{\omega},
                                                                                    & \quad  \textrm{in}\quad Q_T,\\
& Y_k=0,                                                              & \quad  \textrm{on}\quad \Sigma_T,\\ 
& Y_k(\cdot,0)=u_0-y_k(\cdot,0),  Y_{k,t}(\cdot,0)=u_1-y_{k,t}(\cdot,0)                                   & \quad  \textrm{in}\quad \Omega
\end{aligned}
\right.
   \end{equation}
   such that $Y_k(\cdot,T)=-y_k(\cdot,T)$ and $Y_{k,t}(\cdot,T)=-y_{k,t}(\cdot,T)$ in $\Omega$.
This linearization makes appear an operator $K_N$, so that $y_{k+1}=K_N(y_k)$, involving the first derivative of $g$. 
 However, as is well known, such sequence may not converge if the initial guess $(y_0,f_0)$ is not close enough to a zero of $F$ (we refer again to \cite{EFC-AM-2012} exhibiting divergence of the sequence for large data).
  
%
The controllability of nonlinear partial differential equations has attracted a large number of works in the last decades. We refer to the monography \cite{Coron-Book-07} and the references therein.
However, as far as we know, very few are concerned with the approximation of exact controls for nonlinear partial differential equations, so that the construction of convergent control approximations for controllable nonlinear equation remains a challenges. 

Assuming that  $g^{\prime}\in L^\infty_{loc}(\mathbb{R})$ satisfies an asymptotic property similar to \eqref{asymptotic_behavior}, and in addition that there exists one $r$ in $(0,1]$ such that $\sup_{a,b\in \mathbb{R}, a\neq b} \frac{\vert g^\prime(a)-g^\prime(b)\vert}{\vert a-b\vert^r}< \infty$, we construct, for any initial data $(u_0,u_1)\in \boldsymbol{V}$, a strongly convergent sequence $(f_k)_{k\in \mathbb{N}}$ toward a control for (\ref{eq:wave-NL}). Moreover, after a finite number of iterates, the convergence is super linear with a rate equal to $1+r$. This is done  by introducing a quadratic functional which measures how a pair $(y,f) \in Y$ is close to a controlled solution for (\ref{eq:wave-NL}) and then by determining a particular minimizing sequence enjoying the announced property. A natural example of so-called error (or least-squares) functional is given by $\widetilde{E}(y,f):=\frac{1}{2}\Vert \widetilde{F}(y,f)\Vert^2_W$ to be minimized over $Y$. In view of controllability results for~(\ref{eq:wave-NL}), the non-negative functional $\widetilde{E}$ achieves its global minimum equal to zero for any control pair $(y,f)\in Y$ of (\ref{eq:wave-NL}). Inspired by recent works concerning the Navier-Stokes system (see \cite{lemoine-Munch-Pedregal-AMO-20}), we determine, through a specific descent direction, a minimizing sequence $(y_k,f_k)_{k\geq 0}$ converging to a zero of the quadratic functional. 
 
The paper is organized as follows.  Then, in Section \ref{sec:LS}, we define the least-squares functional $E$ and the corresponding optimization problem \eqref{extremal_problem} over the Hilbert $\mathcal{A}$. We show that $E$ is Gateaux-differentiable differentiable over $\mathcal{A}$ and that any critical point $(y,f)$ for $E$ for which $g^\prime(y)$ belongs to $L^\infty(Q_T)$ is also a zero of $E$. This is done by introducing a descent direction $(Y^1,F^1)$ for $E$ at any $(y,f)$ for which $E^\prime(y,f)\cdot (Y^1,F^1)$ is proportional to $E(y,f)$.  Then, assuming that the nonlinear function $g$ satisfies the above conditions, notably that  $\sup_{a,b\in \mathbb{R}, a\neq b} \frac{\vert g^\prime(a)-g^\prime(b)\vert}{\vert a-b\vert^r}< \infty$ for some $r$ in $(0,1]$, we determine a minimizing sequence based on $(Y^1,F^1)$ which converges strongly to a controlled pair for the semilinear wave equation (\ref{eq:wave-NL}). Moreover, we prove that after a finite number of iterates, the convergence enjoys a rate equal to $1+s$. We also emphasize in Section \ref{sec:remarks} that this least-squares approach coincides with the damped Newton method one may use to find a zero of a mapping similar to $\widetilde{F}$ mentioned above. The appendix section \ref{sec:linearizedwave} states some a priori estimates for the linearized wave equation with potential in $L^\infty(Q_T)$ and source term in $L^2(Q_T)$ and emphasize that the operator $K$ is contractant under smallness assumption on $\Vert \hat{g}^\prime\Vert_{L_{loc}^\infty(\mathbb{R})}$

  As far as we know, the method introduced and analyzed in this work is the first one providing an explicit construction of controls for semilinear wave equation. 

 Along the text,  we shall denote by $\Vert \cdot \Vert_{\infty} $ the usual norm in $L^\infty(\mathbb{R})$, $(\cdot,\cdot)_X$ the scalar product of $X$ (if $X$ is a Hilbert space) and by $\langle \cdot, \cdot \rangle_{X,Y}$ the duality product between the spaces $X$ and $Y$.
We shall also denote by $C=C(\Omega,T), C_1=C_1(\Omega,T)$, ..., positive constants only dependent on $\Omega$ and $T$. Last, we shall use the notation $\Vert \cdot\Vert_{2,q_T}$ for $\Vert \cdot\Vert_{L^2(q_T)}$ and  $\Vert \cdot\Vert_p$ for $\Vert \cdot\Vert_{L^p(Q_T)}$, mainly for $p=2$ and $p=\infty$.

\section{The least-squares method and its analysis}\label{sec:LS}

For any $s\in [0,1]$, we define the space 
$$
W_s=\biggl\{g\in C(\R),  \ g^\prime\in L^\infty_{loc}(\mathbb{R}), \sup_{a,b\in \mathbb{R}, a\neq b} \frac{\vert g^\prime(a)-g^\prime(b)\vert}{\vert a-b\vert^s}< \infty\biggr\}.
$$ 
The case $s=0$ reduces to $W_0=\{g\in {\mathcal{C}}(\R), \ g^\prime\in L^\infty_{loc}(\mathbb{R})\}$ while the case $s=1$ corresponds to $
W_1=\{g\in {\mathcal{C}}(\R), \ g^\prime\in L^\infty_{loc}(\mathbb{R}), g^{\prime\prime}\in L^\infty(\mathbb{R})\}$.

\subsection{The least-squares method}

We assume that $g$ belongs to $W_0$ and introduce the vector space $\mathcal{A}_0$ 
\begin{equation}
\nonumber
\begin{aligned}
\mathcal{A}_0=
\biggl\{(y,f): \,& y\in L^2(Q_T), (y(\cdot,0),y_t(\cdot,0))\in \boldsymbol{V}, f\in L^2(q_T), y_{tt} - y_{xx}\in L^2(Q_T),  \\
& (y(\cdot,0),y_t(\cdot,0))=(0,0),\, (y(\cdot,T),y_t(\cdot,T))=(0,0)\, \textrm{in}\,\, \Omega, \, y=0 \,\textrm{on}\, \Sigma_T\biggr\}.
\end{aligned}
\end{equation}
Endowed with the scalar product
\begin{equation}\nonumber
\begin{aligned}
((y,f),(\overline{y},\overline{f}))_{\mathcal{A}_0}= &(y,\overline{y})_2+ ((y(\cdot,0),y_t(\cdot,0)),(\overline{y}(\cdot,0),\overline{y}_t(\cdot,0)))_{\boldsymbol{V}}\\
&+(y_{tt} - y_{xx},\overline{y}_{tt} - \overline{y}_{xx})_2+ (f,\overline{f})_{2,q_T}
\end{aligned}
\end{equation}
$\mathcal{A}_0$ is an Hilbert space. We shall note $\Vert (y,f)\Vert_{\mathcal{A}_0}:=\sqrt{((y,f),(y,f))_{\mathcal{A}_0}}$. We also consider the affine (convex) space $\mathcal{A}$  
\begin{equation}
\nonumber
\begin{aligned}
\mathcal{A}=
\biggl\{(y,f): \,& y\in L^2(Q_T), (y(\cdot,0),y_t(\cdot,0))\in \boldsymbol{V}, f\in L^2(q_T), y_{tt} - y_{xx}\in L^2(Q_T),  \\
& (y(\cdot,0),y_t(\cdot,0))=(u_0,u_1),\, (y(\cdot,T),y_t(\cdot,T))=(z_0,z_1)\, \textrm{in}\,\, \Omega, \, y=0 \,\textrm{on}\, \Sigma_T\biggr\}.
\end{aligned}
\end{equation}
Observe also that we can write $\mathcal{A}=(\overline{y},\overline{f})+\mathcal{A}_0$ for any element $(\overline{y},\overline{f})\in \mathcal{A}$.

For any fixed $(\overline{y},\overline{f})\in \mathcal{A}$, we consider the following extremal problem :
\begin{equation}
\label{extremal_problem}
\inf_{(y,f)\in \mathcal{A}_0} E(\overline{y}+y,f+\overline{f})
\end{equation}
where $E:\mathcal{A}\to \mathbb{R}$ is defined as follows 
\begin{equation}\nonumber
E(y,f):=\frac{1}{2}\big\Vert y_{tt}-y_{xx} + g(y)-f\,1_{\omega}\big\Vert^2_{L^2(Q_T)}
\end{equation}
justifying the least-squares terminology we have used. 

Remark that the functional $E$ is well-defined in $\mathcal{A}$. Precisely, \textit{a priori} estimate for the linear wave equation reads as  
$$
\Vert (y,y_t)\Vert^2_{L^{\infty}(0,T; \boldsymbol{V})} \leq C \biggl( \Vert y_{tt}-y_{xx}\Vert^2_{L^2(Q_T)}+\Vert y_0,y_1\Vert^2_{\boldsymbol{V}}\biggr)
$$
for any $y$ such that $(y,f)\in \mathcal{A}$  and implies in this one dimensional setting that $y\in L^\infty(Q_T)$. Since $g\in W_0$, $\vert g(s)-g(0)\vert \leq \Vert g^{\prime}\Vert_{\infty} \vert s\vert$ for all $s\in \mathbb{R}$ so that $\Vert g(y)\Vert_2\leq \vert g(0)\vert \sqrt{\vert Q_T\vert}+ \Vert g^{\prime}\Vert_{\infty,loc} \Vert y\Vert_2$.

Within the hypotheses of Theorem \ref{zuazuaTH}, the infimum of the functional of $E$ is zero and is reached by at least one pair $(y,f)\in \mathcal{A}$, solution of  (\ref{eq:wave-NL}) and satisfying $(y(\cdot,T),y_t(\cdot,T))=(z_0,z_1)$. Conversely, any pair $(y,f)\in \mathcal{A}$ for which $E(y,f)$ vanishes is solution of (\ref{eq:wave-NL}). In this sense, the functional  $E$ is a so-called error functional which measures the deviation of $(y,f)$ from being a solution of the underlying nonlinear equation.  A practical way of taking a functional to its minimum is through some clever use of descent  directions, i.e the use of its derivative. In doing so, the presence of local minima is always something that  may dramatically spoil the whole scheme. The unique structural property that discards this possibility is the strict convexity of the functional $E$. However, for nonlinear equation like (\ref{eq:wave-NL}), one cannot expect this property to hold for the functional $E$. Nevertheless, we insist in that one may construct a particular minimizing sequence which cannot converge except to a global minimizer leading $E$ down to zero. 

In order to construct such minimizing sequence, we look, for any $(y,f)\in \mathcal{A}$, for a pair $(Y^1,F^1)\in \mathcal{A}_0$ solution of the following formulation 

\begin{equation}
\label{wave-Y1}
\left\{
\begin{aligned}
& Y^1_{tt} - Y^1_{xx} +  g^{\prime}(y)\cdot Y^1 = F^1 1_{\omega}+\big(y_{tt}-y_{xx} + g(y)-f\, 1_{\omega}\big), &\textrm{in}\,\, Q_T,\\
& Y^1=0, &\textrm{on}\,\, \Sigma_T, \\
& (Y^1(\cdot,0),Y^1_t(\cdot,0))=(0,0), & \textrm{in}\,\, \Omega.
\end{aligned}
\right.
\end{equation}
Remark that $(Y^1,F^1)$ belongs to $\mathcal{A}_0$ if and only if  $F^1$
is a null control for $Y^1$. Among the controls of this linear equation, we select the control of minimal $L^2(q_T)$ norm. We have the following property. 

\begin{lemma}\label{estimateC1C2} Assume $g\in W_0$. Assume that $T>2\max(l_1,1-l_2)$. Let any $(y,f)\in \mathcal{A}$. There exists a pair $(Y^1,F^1)\in \mathcal{A}_0$ solution of (\ref{wave-Y1}). Moreover, the pair $(Y^1,F^1)$ for which $\Vert F^1\Vert_{2,q_T}$ is minimal satisfies the following estimates :  

\begin{equation} \label{estimateF1Y1}
\Vert (Y^1,Y^1_t)\Vert_{L^{\infty}(0,T;\boldsymbol{V})}+ \Vert F^1\Vert_{2,q_T} \leq C_1  e^{C_2\sqrt{\Vert g^{\prime}(y)\Vert_{\infty}} }\sqrt{E(y,f)},
\end{equation}

and 

\begin{equation} \label{estimateF1Y1_A0}
\Vert (Y^1,F^1)\Vert_{\mathcal{A}_0} \leq C_1  e^{C_2\sqrt{\Vert g^{\prime}(y)\Vert_{\infty}} }\sqrt{E(y,f)}.
\end{equation}
\end{lemma}
\textsc{Proof-} The first estimate is a consequence of Proposition \ref{controllability_result} (see the appendix) using the equality $\Vert \partial_{tt} y-\Delta y +g(y)-f\,1_{\omega}\Vert_2=\sqrt{2E(y,f)}$. The second one follows from 
\begin{equation}
\nonumber
\begin{aligned}
\Vert (Y^1,F^1)\Vert_{\mathcal{A}_0} & \leq  \Vert Y^1_{tt}-Y^1_{xx}\Vert_2 + \Vert Y^1\Vert_2+\Vert F^1\Vert_{2,q_T}+\Vert Y^1(\cdot,0),Y^1_t(\cdot,0)\Vert_{\boldsymbol{V}} \\
& \leq (1+\Vert g^{\prime}(y)\Vert_\infty)\Vert Y^1\Vert_2 + 2\Vert F^1\Vert_{2,q_T}+\sqrt{E(y,f)}\\
& \leq C_1 (1+\Vert g^{\prime}(y)\Vert_\infty)e^{C_2\sqrt{\Vert g^{\prime}(y)\Vert_{\infty}}}\sqrt{E(y,f)}\\
& \leq C_1 e^{(1+C_2)\sqrt{\Vert g^{\prime}(y)\Vert_{\infty}}}\sqrt{E(y,f)}.
\end{aligned}
\end{equation}
using that $(1+s)\leq e^{2\sqrt{s}}$ for all $s\geq 0$.
$\hfill\Box$

In particular, this implies that 
\begin{equation}\nonumber
\Vert Y^1\Vert_{L^{\infty}(Q_T)} \leq C_1  e^{C_2\sqrt{\Vert g^{\prime}(y)\Vert_{\infty}} }\sqrt{E(y,f)}.
\end{equation}

The interest of the pair $(Y^1,F^1)\in \mathcal{A}_0$ lies in the following result.
\begin{lemma}\label{differentiabiliteE}
Assume that $g\in W_0$ and $T>2\max(l_1,1-l_2)$. Let $(y,f)\in \mathcal{A}$  and let $(Y^1,F^1)\in \mathcal{A}_0$ be a solution of (\ref{wave-Y1}). Then the derivative of $E$ at the point $(y,f)\in \mathcal{A}$ along the direction $(Y^1,F^1)$ defined by $E^{\prime}(y,f)\cdot (Y^1,F^1):= \lim_{\lambda\to 0, \lambda\neq 0} \frac{E((y,f)+\lambda (Y^1,F^1))-E(y,f)}{\lambda}$ satisfies 
\begin{equation}\label{estimateEEprime}
E^{\prime}(y,f)\cdot (Y^1,F^1)=2E(y,f).
\end{equation}
\end{lemma} 
\textsc{Proof-} We preliminary check that for all $(Y,F)\in \mathcal{A}_0$, the functional $E$ is differentiable at the point $(y,f)\in \mathcal{A}$ along the direction $(Y,F)\in \mathcal{A}_0$. For any $\lambda\in \mathbb{R}$, simple computations lead to the equality 

$$
\begin{aligned}
E(y+\lambda Y,f+\lambda F) =E(y,f)+ \lambda E^{\prime}(y,f)\cdot (Y,F) + h((y,f),\lambda (Y,F))
\end{aligned}
$$
with 
\begin{equation}\label{Efirst}
E^{\prime}(y,f)\cdot (Y,F):=\biggl(y_{tt}-y_{xx}+ g(y)-f\, 1_\omega,  Y_{tt}-Y_{xx}+ g^\prime(y)Y-F\, 1_\omega\biggr)_2
\end{equation}
and
$$
\begin{aligned}
h((y,f),\lambda (Y,F)):=& \frac{\lambda^2}{2}\biggl(Y_{tt}-Y_{xx}+ g^\prime(y)Y-F\, 1_\omega,Y_{tt}-Y_{xx}+ g^\prime(y)Y-F\, 1_\omega\biggl)_2\\
& +\lambda \biggl(Y_{tt}-Y_{xx}+ g^\prime(y)Y-F\, 1_\omega,l(y,\lambda Y)\biggl)_2\\
& +\biggl(y_{tt}-y_{xx}+ g(y)-f\, 1_\omega,l(y,\lambda Y)\biggr)+ \frac{1}{2}(l(y,\lambda Y),l(y,\lambda Y))
\end{aligned}
$$ 
where $l(y,\lambda Y):=g(y+\lambda Y)-g(y)-\lambda g^{\prime}(y)Y$.
The application $(Y,F)\to E^{\prime}(y,f)\cdot (Y,F)$ is linear and continuous from $\mathcal{A}_0$ to $\mathbb{R}$ as it satisfies
\begin{equation}
\label{useful_estimate}
\begin{aligned}
\vert E^{\prime}(y,f)\cdot (Y,F)\vert &\leq  \Vert y_{tt}-y_{xx}+ g(y)-f\, 1_\omega\Vert_2 \Vert Y_{tt}-Y_{xx}+ g^\prime(y)Y-F\, 1_\omega\Vert_2\\
& \leq \sqrt{2E(y,f)} \biggl(\Vert (Y_{tt}-Y_{xx})\Vert_2 + \Vert g^\prime(y)\Vert_{\infty}\Vert Y\Vert_2 + \Vert F\Vert_{2,q_T}\biggr)\\
& \leq \sqrt{2E(y,f)} \max\big(1,\Vert g^\prime(y)\Vert_\infty\big) \Vert (Y,F)\Vert_{\mathcal{A}_0}.
\end{aligned}
\end{equation}
Similarly, for all $\lambda\in \mathbb{R}^\star$, 
$$
\begin{aligned}
\biggl\vert \frac{1}{\lambda} h((y,f),\lambda (Y,F))&\biggr\vert  \leq \frac{\lambda}{2} \Vert Y_{tt}-Y_{xx}+ g^\prime(y)Y-F\, 1_\omega\Vert^2_2 \\
& +\biggl(\lambda\Vert Y_{tt}-Y_{xx}+ g^\prime(y)Y-F\, 1_\omega\Vert_2 +\sqrt{2E(y,f)}+\frac{1}{2}\Vert l(y,\lambda Y)\Vert_2\biggl)\frac{1}{\lambda}\Vert l(y,\lambda Y)\Vert_2.
\end{aligned}
$$
Since $g^\prime\in L^\infty_{loc}(\mathbb{R})$ and $y\in L^{\infty}(Q_T)$, we have
$$
\biggl\vert \frac{1}{\lambda}l(y,\lambda Y)\biggl\vert = \biggl\vert \frac{g(y+\lambda Y)-g(y)}{\lambda}-g^{\prime}(y)Y\biggr\vert\leq 2\Vert g^{\prime}(y)\Vert_{L^\infty(Q_T)} \vert Y\vert, \quad a.e. \,\, \textrm{in}\,\,Q_T
$$
and that $\biggl\vert \frac{1}{\lambda}l(y,\lambda Y)\biggl\vert = \biggl\vert \frac{g(y+\lambda Y)-g(y)}{\lambda}-g^{\prime}(y)Y\biggr\vert\to 0$ as $\lambda\to 0$, a.e. in $Q_T$. From the Lebesgue's Theorem, it follows that $\vert \frac{1}{\lambda}\vert \Vert l(y,\lambda Y)\Vert_2\to 0$ as $\lambda\to 0$ and then that $\vert h((y,f),\lambda (Y,F))\vert=o(\lambda)$. 
We deduce that the functional $E$ is differentiable at the point $(y,f)\in \mathcal{A}$ along the direction $(Y,F)\in \mathcal{A}_0$.

Eventually, the equality (\ref{estimateEEprime}) follows from the definition of the pair $(Y^1,F^1)$ given in (\ref{wave-Y1}). $\hfill\Box$


\vskip 0.25cm
Remark that from the equality \eqref{Efirst}, the derivative $E^{\prime}(y,f)$ is independent of $(Y,F)$. We can then define the norm $\Vert E^{\prime}(y,f)\Vert_{(\mathcal{A}_0)^{\prime}}:= \sup_{(Y,F)\in \mathcal{A}_0} \frac{E^{\prime}(y,f)\cdot (Y,F)}{\Vert (Y,F)\Vert_{\mathcal{A}_0}}$ associated to $(\mathcal{A}_0)^{\prime}$, the set of the linear and continuous applications from $\mathcal{A}_0$ to $\mathbb{R}$

\vskip 0.25cm

Combining the equality \eqref{estimateEEprime} and the inequality \eqref{estimateF1Y1}, we deduce the following estimate of $E(y,f)$ in term of the norm of $E^\prime(y,f)$. 
\begin{prop}Assume $g\in W_0$ and $T>2\max(l_1,1-l_2)$. For any $(y,f)\in\mathcal{A}$, the following inequalities hold true:
\begin{equation}\label{ineq_E_Eprime}
\frac{1}{\sqrt{2}\max\big(1,\Vert g^\prime(y)\Vert_\infty\big)}\Vert E^{\prime}(y,f)\Vert_{\mathcal{A}_0^\prime}\leq \sqrt{E(y,f)}\leq \frac{1}{\sqrt{2}}C_1 e^{C_2\sqrt{\Vert g^{\prime}(y)\Vert_\infty}} \Vert E^{\prime}(y,f)\Vert_{\mathcal{A}_0^\prime}
\end{equation}
where $C_1$ and $C_2$ are the positive constants from Proposition \ref{estimateC1C2}.
\end{prop}
\textsc{Proof-} \eqref{estimateEEprime} rewrites
$E(y,f)=\frac12 E^{\prime}(y,f)\cdot (Y^1,F^1)$
where $(Y^1,F^1)\in \mathcal{A}_0$ is solution of (\ref{wave-Y1})  and therefore, with \eqref{estimateF1Y1_A0}
$$E(y,f)\le\frac12 \|E'(y,f)\|_{ \mathcal{A}_0'} \|(Y^1,F^1)\|_{ \mathcal{A}_0}\le \frac12  C_1 e^{C_2\sqrt{\Vert g^{\prime}(y)\Vert_\infty}} \|E'(y,f)\|_{ \mathcal{A}_0'}\sqrt{E(y,f)}.
$$
On the other hand, for all $(Y,F)\in\mathcal{A}_0$, the inequality \eqref{useful_estimate}, i.e.
$$
\vert E^{\prime}(y,f)\cdot (Y,F)\vert \leq  \sqrt{2E(y,f)} \max\biggl(1,\Vert  g^\prime(y)\Vert_{\infty}\biggr) \Vert (Y,F)\Vert_{\mathcal{A}_0}
$$
leads to the left inequality.$\hfill\Box$

In particular, any \textit{critical} point $(y,f)\in \mathcal{A}$ for $E$ (i.e. for which $E^\prime(y,f)$ vanishes) such that $\Vert g^\prime(y)\Vert_{L^\infty(Q_T)}<\infty$ is a zero for $E$, a pair solution of the controllability problem. 
In other words, any sequence $(y_k,f_k)_{k>0}$ satisfying $\Vert E^\prime(y_k,f_k)\Vert_{\mathcal{A}_0^\prime}\to 0$ as $k\to \infty$ and for which $\Vert g^{\prime}(y_k)\Vert_{\infty}$
is uniformly bounded is such that $E(y_k,f_k)\to 0$ as $k\to \infty$. We insist that this property does not imply the convexity of the functional $E$ (and \textit{a fortiori} the strict convexity of $E$, which actually does not hold here in view of the multiple zeros for $E$) but show that a minimizing sequence for $E$ can not be stuck in a local minimum. 

On the other hand, the left inequality indicates the functional $E$ is flat around its zero set. As a consequence, gradient based minimizing sequences may achieve a very low rate of convergence (we refer to \cite{AM-PP-2014}
and also \cite{lemoine-Munch-Pedregal-AMO-20} devoted to the Navier-Stokes equation where this phenomenon is observed).

\subsection{A strongly convergent minimizing sequence for $E$}

\ 
 We now examine the convergence of an appropriate sequence $(y_k,f_k)\in \mathcal{A}$. 
 In this respect, we observe that equality (\ref{estimateEEprime}) shows that $-(Y^1,F^1)$ given by the solution of (\ref{wave-Y1}) is a descent direction for the functional $E$. Therefore, we can define at least formally, for any $m\geq 1$,  a minimizing sequence $(y_k,f_k)_{k>0}$  as follows: 
\begin{equation}
\label{algo_LS_Y}
\left\{
\begin{aligned}
&(y_0,f_0) \in \mathcal{A}, \\
&(y_{k+1},f_{k+1})=(y_k,f_k)-\lambda_k (Y^1_k,F_k^1), \quad k>0, \\
& \lambda_k= \textrm{argmin}_{\lambda\in (0,m]} E\big((y_k,f_k)-\lambda (Y^1_k,F_k^1)\big),    
\end{aligned}
\right.
\end{equation}
where $(Y^1_k,F_k^1)\in \mathcal{A}_0$ is such that $F^1_k$ is the null control of minimal $L^2(q_T)$-norm for $Y^1_k$, solution of 
\begin{equation}
\label{wave-Y1k}
\left\{
\begin{aligned}
& Y_{k,tt}^1 - Y^1_{k,xx} +  g^{\prime}(y_k)\cdot Y^1_k = F^1_k 1_{\omega}+ (y_{k,tt}-y_{k,xx}+g(y_k)-f_k 1_\omega), & \textrm{in}\,\, Q_T,\\
& Y_k^1=0,  & \textrm{on}\,\, \Sigma_T, \\
& (Y_k^1(\cdot,0),Y_{k,t}^1(\cdot,0))=(0,0), & \textrm{in}\,\, \Omega.
\end{aligned}
\right.
\end{equation}

We prove in this section the strong convergence of the sequence $(y_k,f_k)_{k\in\mathbb{N}}$ toward a controlled pair for \ref{eq:wave-NL}, first in the case $g^\prime\in L^\infty(\mathbb{R})$ in Theorem \ref{2.12} and then in the case $g^\prime\in L^\infty(\mathbb{R})$ and satisfying a growth condition at infinity in Theorem \ref{Th2.12}

We first perform the analysis assuming the non linear function $g$ in $W_1$, notably that $g^{\prime\prime}\in L^{\infty}(\mathbb{R})$ (the derivatives here are in the sense of distribution). We first prove the following lemma.

\begin{lemma}\label{lemma_estimateE} Assume that $g\in W_1$ and that $T>2\max(l_1,1-l_2)$. For any $(y,f)\in \mathcal{A}$, let $(Y^1,F^1)\in \mathcal{A}_0$ be defined by \eqref{wave-Y1}. For any $\lambda\in \mathbb{R}$ and $k\in \mathbb{N}$, the following estimate holds, for some $C>0$
\begin{equation}\label{estimateE}
\left\{
\begin{aligned}
&E\big((y,f)-\lambda (Y^1,F^1)\big)\leq   E(y,f)\biggl(\vert 1-\lambda\vert +\lambda^2\,C(y) \sqrt{E(y,f)}\biggr)^2,\\
& C(y):= \frac{C}{2\sqrt{2}} \Vert g^{\prime\prime}\Vert_{\infty}  \biggl(C_1  e^{C_2\sqrt{\Vert g^{\prime}(y)\Vert_\infty}}\biggr)^2.
\end{aligned}
\right.
\end{equation}
\end{lemma}
\textsc{Proof-} With $g\in W_1$, we write that 
\begin{equation}
\label{majoration-l}
\vert l(y_,-\lambda Y^1)\vert = \vert g(y-\lambda Y^1)-g(y)+\lambda g^{\prime}(y) Y^1\vert \leq \frac{\lambda^2}{2}\Vert g^{\prime\prime}\Vert_\infty (Y^1)^2
\end{equation}
and obtain that
\begin{equation}
\label{expansionE}
\begin{aligned}
2 &E\big((y,f)-\lambda (Y^1,F^1)\big) \\
&=\biggl\Vert \big(y_{tt}-y_{xx}+g(y)-f\,1_\omega\big)-\lambda \big(Y^1_{tt}-Y^1_{xx}+g^\prime(y)Y^1-F\,1_\omega\big)+l(y,-\lambda Y^1)\biggr\Vert^2_2\\
& =\biggl\Vert (1-\lambda)\big(y_{tt}-y_{xx}+g(y)-f\,1_\omega\big)+l(y,-\lambda Y^1)\biggr\Vert^2_2\\
& \leq \biggl( \big\Vert (1-\lambda)\big(y_{tt}-y_{xx}+g(y)-f\,1_\omega\big)\big\Vert_2 +\big\Vert l(y,-\lambda Y^1)\big\Vert_2\biggr)^2\\
&\leq 2 \biggl( \vert 1-\lambda\vert \sqrt{E(y,f)}+\frac{\lambda^2}{2\sqrt{2}}\Vert g^{\prime\prime}\Vert_\infty \Vert(Y^1)^2\Vert_2\biggr)^2.
\end{aligned}
\end{equation}
But, in view of  \eqref{estimateF1Y1}, we have 
\begin{equation}
\label{estim1}
\begin{aligned}
\Vert(Y^1)^2\Vert_2 & \leq \Vert Y^1\Vert_{L^\infty(Q_T)}\Vert Y^1\Vert_{L^2(Q_T)} \leq C \Vert Y^1\Vert^2_{L^{\infty}(0,T; H_0^1(\Omega))} \\
& \leq C \biggl(C_1e^{C_2\sqrt{\Vert g^{\prime}(y)\Vert_{\infty}}}\biggr)^2 E(y,f) 
\end{aligned}
\end{equation}
for some constant $C=C(\Omega,T)$ from which we get \eqref{estimateE}. $\hfill\Box$

The previous result still holds if we assume only that $g\in W_s$ for some $s\in (0,1)$. For any $g\in W_s$, we introduce the notation 
$\|g^\prime\|_{\widetilde{W}^{s,\infty}(\R)}:=\sup_{a,b\in \mathbb{R}, a\neq b} \frac{\vert g^\prime(a)-g^\prime(b)\vert}{\vert a-b\vert^s}$. We have the following result. 

\begin{lemma}\label{lemma_estimateEs} Assume that $g\in W_s$ for some $s\in (0,1)$ and that $T>2\max(l_1,1-l_2)$. For any $(y,f)\in \mathcal{A}$, let $(Y^1,F^1)\in \mathcal{A}_0$ be defined by \eqref{wave-Y1}. For any $\lambda\in \mathbb{R}$ and $k\in \mathbb{N}$, the following estimate holds
\begin{equation}\label{estimateEs}
\left\{
\begin{aligned}
&E\big((y,f)-\lambda (Y^1,F^1)\big)\leq   E(y,f)\biggl(\vert 1-\lambda\vert +\lambda^{1+s}\,C(y) E(y,f)^{s/2}\biggr)^2,\\
& C(y):= \frac{C}{2\sqrt{2}} \|g^\prime\|_{\widetilde{W}^{s,\infty}(\R)}  \biggl(C_1  e^{C_2\sqrt{\Vert g^{\prime}(y)\Vert_\infty}}\biggr)^{1+s}.
\end{aligned}
\right.
\end{equation}
\end{lemma}
\textsc{Proof-} For any $(x,y)\in \R^2$ and $\lambda\in \R$, we write $g(x+\lambda y)-g(x)=\int_0^\lambda y g'(x+\xi y)d\xi$ leading to 
$$
\begin{aligned}
|g(x+\lambda y)-g(x)-\lambda g'(x)y|
&\le \int_0^\lambda |y| |g'(x+\xi y)-g'(x)|d\xi\\
&\le \int_0^\lambda |y|^{1+s}|\xi|^s\frac{|g'(x+\xi y)-g'(x)|}{|\xi y|^s}d\xi\\
&\le \|g^\prime\|_{\widetilde{W}^{s,\infty}(\R)}|y|^{1+s}\frac{\lambda^{1+s}}{1+s}.
\end{aligned}
$$
It follows that 
$$|l(y,-\lambda Y^1)|=|g(y-\lambda Y^1)-g(y)+\lambda g^{\prime}(y) Y^1|\le  \|g^\prime\|_{\widetilde{W}^{s,\infty}(\R)}\frac{\lambda^{1+s}}{1+s}|Y^1|^{1+s}
$$
and 
$$\begin{aligned}
\bigl \Vert l(y,\lambda Y^1)\bigr\Vert_2
&\le  \|g^\prime\|_{\widetilde{W}^{s,\infty}(\R)}\frac{\lambda^{1+s}}{1+s}\bigl \Vert |Y^1|^{1+s}\bigr\Vert_{L^2(0,T;L^{2}(\Omega))}.
\end{aligned}
 $$
 But 
 $$
 \bigl \Vert |Y^1|^{1+s}\bigr\Vert^2_{L^2(0,T;L^{2}(\Omega))}=\int_{Q_T}(Y^1)^{2(s+1)} \, dx\, dt =\Vert Y^1\Vert_{L^{2(s+1)}(Q_T)}^{2(s+1)}\leq C \Vert Y^1\Vert_{L^{\infty}(Q_T)}^{2(s+1)}
$$
leading to 
 $$
 \bigl \Vert |Y^1|^{1+s}\bigr\Vert_{L^2(0,T;L^{2}(\Omega))}\leq C \biggl(C_1e^{C_2\sqrt{\Vert g^{\prime}(y)\Vert_{\infty}}}\biggr)^{1+s} E(y,f)^{\frac{1+s}{2}}
 $$
 and to the result. $\hfill\Box$
 
 Proceeding as in \cite{ lemoinemunch_time}, we are now in position to prove the strong convergence result for the sequences $(E(y_k,f_k))_{(k\geq 0)}$ and $(y_k,f_k)_{(k\geq 0)}$ for the norm $\Vert \cdot\Vert_\A$.
In order to fix notations and arguments, we first start by making the stronger assumption that $g^\prime\in L^\infty(\mathbb{R})$.
 
 \subsubsection{Convergence in the case $g\in W_s$, $s\in (0,1]$ and $g^\prime\in L^\infty(\mathbb{R})$}

\begin{prop}\label{geomquadratic}
Assume $g\in W_s$ for some $s\in (0,1]$ and that $g^\prime\in L^\infty(\mathbb{R})$.
Let $(y_k,f_k)_{k>0}$ be the sequence of $\mathcal{A}$ defined in \eqref{algo_LS_Y}. Then $E(y_k,f_k)\to 0$ as $k\to\infty$. Moreover, there exists a $k_0\in \mathbb{N}$ such that the sequence $(E(y_k,f_k))_{k>k_0}$ decays with a rate equal to $s+1$. 
\end{prop}

\textsc{Proof-} We make the proof in the case $s\in (0,1)$. The proof in the case $s=1$ is very similar (see next section). Since $g^{\prime}\in L^\infty(\mathbb{R})$, the constant
$C(y_k)$ in \eqref{estimateEs} is uniformly bounded w.r.t. $k$ so that, for all $k>0$ $C(y_k)\leq c$ for some $c>0$.
 For any $(y_k,f_k)\in \mathcal{A}$, let us then denote the real function $p_k$ by 
$$p_k(\lambda):=\vert 1-\lambda\vert +\lambda^{1+s} c\, E(y_k,f_k)^{s/2}$$
 for all $\lambda\in [0,m]$. Lemma \ref{lemma_estimateE} with $(y,f)=(y_k,f_k)$ then allows to write that 
 \begin{equation}\label{estek}
\sqrt{E(y_{k+1},f_{k+1})}=\min_{\lambda\in [0,m]} \sqrt{E((y_k,f_k)-\lambda (Y^1_k,F^1_k))}\leq \min_{\lambda\in [0,m]} p_k(\lambda) \sqrt{E(y_k,f_k)}=p_k(\widetilde{\lambda_k}) \sqrt{E(y_k,f_k)}.
\end{equation}
We then easily check that the optimal $\widetilde{\lambda_k}$ is given by 
\begin{equation}\nonumber
\widetilde{\lambda_k}:=\left\{
\begin{aligned}
& \frac{1}{(1+s)^{1/s} c^{1/s}\sqrt{E(y_k,f_k)}}, & \textrm{if}\quad (1+s)^{1/s} c^{1/s}\sqrt{E(y_k,f_k)}\geq 1,\\
& 1, & \textrm{if}\quad (1+s)^{1/s} c^{1/s}\sqrt{E(y_k,f_k)}< 1
\end{aligned}
\right.
\end{equation}
leading to 
\begin{equation} \label{ptildek}
p_k(\widetilde{\lambda_k}):=\left\{
\begin{aligned}
& 1- \frac{s}{(1+s)^{\frac1s+1}}\frac{1}{c^{1/s}\sqrt{E(y_k,f_k)}}, & \textrm{if}\quad (1+s)^{1/s} c^{1/s}\sqrt{E(y_k,f_k)}\geq 1,\\
& c\, E(y_k,f_k)^{s/2}, & \textrm{if}\quad (1+s)^{1/s} c^{1/s}\sqrt{E(y_k,f_k)}< 1.
\end{aligned}
\right.
\end{equation}

If $(1+s)^{1/s} c^{1/s}\sqrt{E(y_0,f_0)}< 1$, then  $c^{1/s}\sqrt{E(y_0,f_0)}< 1$ (and thus $c^{1/s}\sqrt{E(y_k,f_k)}<1$ for all $k\in\mathbb{N}$) then \eqref{estek} implies that 
$$
c^{1/s}\sqrt{E(y_{k+1},f_{k+1})}\le \big(c^{1/s}\sqrt{E(y_k,f_k)}\big)^{1+s}.
$$
It follows that $c^{1/s}\sqrt{E(y_k,f_k)}\to 0$ as $k\to \infty$ with a rate equal to $1+s$.

If $(1+s)^{1/s} c^{1/s}\sqrt{E(y_0,f_0)}\geq 1$ then we check that the set $I:=\{k\in \mathbb{N},\ (1+s)^{1/s} c^{1/s}\sqrt{E(y_k,f_k)}\geq 1\}$ is a finite subset of $\mathbb{N}$;  indeed, for all $k\in I$, \eqref{estek} implies that 
\begin{equation}\label{decayEunquarts}
c^{1/s}\sqrt{E(y_{k+1},f_{k+1}) } \le \Big(1- \frac{s}{(1+s)^{\frac1s+1}}\frac{1}{c^{1/s}\sqrt{E(y_k,f_k)}}\Big)c^{1/s}\sqrt{E(y_k,f_k)}=c^{1/s}\sqrt{E(y_k,f_k) }-\frac{s}{(1+s)^{\frac1s+1}}
\end{equation}
and the strict decrease of the sequence $(c^{1/s} \sqrt{E(y_k,f_k)})_{k\in I}$. Thus there exists 
$k_0\in\mathbb{N}$ such that for all $k\geq k_0$, $(1+s)^{1/s} c^{1/s}\sqrt{E(y_k,f_k)}<1$, that is  $I$ is a finite subset of $\mathbb{N}$.  
Arguing as in the first case, it follows that $\sqrt{E(y_k,f_k)}\to 0$ as $k\to \infty$.

It follows in particular, in view of (\ref{ptildek}) that the sequence $(p_k(\widetilde{\lambda_k}))_{k\in\N}$ decreases as well.  $\hfill\Box$

\begin{remark}
The number of iterates $k_0$ necessary to reach a super-linear regime depends on the value of $E(y_0,f_0)$,  $\Vert g^\prime\Vert_{L^\infty(\mathbb{R})}$ and $\|g^\prime\|_{\widetilde{W}^{s,\infty}(\R)}$. For instance, with $s=1$, writing  from \eqref{decayEunquarts} that $c\sqrt{E(y_k,f_k)}\leq c\sqrt{E(y_0,f_0)}-\frac{k}{4}$ for all $k$ such that $c\sqrt{E(y_k,f_k)}\geq 1$, we obtain that 
$$
k_0\leq  \biggl\lfloor 4(K\sqrt{E(y_0,f_0)}-1)+1 \biggr\rfloor, \quad \textrm{with}\quad K=\mathcal{O}(C\Vert g^{\prime\prime}\Vert_{\infty} e^{C\sqrt{\Vert g^{\prime}\Vert_{\infty}}})
$$
where $\lfloor x\rfloor$ denotes the integer part of $x\in \mathbb{R}^+$. 
\end{remark}

We also have the following convergence of the optimal sequence $(\lambda_k)_{k>0}$. 

\begin{lemma}\label{lambda-k-go-1}
Assume that $g\in W_s$ for some $s\in (0,1]$ and that $g^\prime\in L^\infty(\mathbb{R})$. The sequence $(\lambda_k)_{k>0}$ defined in (\ref{algo_LS_Y}) converges to $1$ as $k\to \infty$.
\end{lemma}
\textsc{Proof-} Take $s=1$. In view of \eqref{expansionE}, we have, as long as  $E(y_k,f_k)>0$, since $\lambda_k\in[0,m]$ 
$$\begin{aligned}
(1-\lambda_k)^2
&=\frac{E(y_{k+1},f_{k+1})}{E(y_{k},f_k)}-2(1-\lambda_k)\frac{\big(y_{k,tt}+y_{k,xx}+g(y_k)-f_k\,1_\omega,l(y_k,\lambda_k Y_k^1) \big)_2}{E(y_{k},f_k)}\\
& \hspace{3cm}- \frac{\bigl \Vert l(y_k,\lambda_k Y_k^1)\bigr\Vert^2_{2}}{2E(y_{k})}\\
&\le \frac{E(y_{k+1},f_{k+1})}{E(y_{k},f_k)}-2(1-\lambda_k)\frac{\big(y_{k,tt}+y_{k,xx}+g(y_k)-f_k\,1_\omega,l(y_k,\lambda_k Y_k^1) \big)_2}{E(y_{k},f_k)}\\
&\le \frac{E(y_{k+1},f_{k+1})}{E(y_{k},f_k)}+2\sqrt{2}m\frac{\sqrt{E(y_k,f_k)}\|l(y_k,\lambda_k Y_k^1) \|_{2}}{E(y_{k},f_k)}\\
&\le \frac{E(y_{k+1},f_{k+1})}{E(y_{k},f_k)}+2\sqrt{2}m\frac{\|\l(y_k,\lambda_k Y_k^1) \|_2}{\sqrt{E(y_k,f_k)}}.
\end{aligned}
$$
But, from \eqref{majoration-l} and  \eqref{estim1} 
$$
\|l(y_k,\lambda_k Y_k^1) \|_{2}
\le \frac{\lambda_k^2}{2\sqrt{2}}\Vert g^{\prime\prime}\Vert_\infty \Vert (Y^1_k)^2\Vert_{2}\le m^2 C(y_k)E(y_k,f_k)
$$
and thus 
$$
(1-\lambda_k)^2\le \frac{E(y_{k+1},f_{k+1})}{E(y_{k},f_k)}+2\sqrt{2}m^3 c \sqrt{E(y_k,f_k)}.
$$
Consequently, since $E(y_{k},f_k)\to 0$ and $\frac{E(y_{k+1},f_{k+1})}{E(y_{k},f_k)}\to 0$, we deduce that $(1-\lambda_k)^2\to 0$. $\hfill\Box$

We are now in position to prove the following convergence result.

 \begin{theorem}\label{2.12}
Assume that $g\in W_s$ for some $s\in (0,1]$ and that $g^\prime\in L^\infty(\mathbb{R})$. Let $(y_k,f_k)_{k\in\mathbb{N}}$ be the sequence defined by (\ref{algo_LS_Y}). Then, $(y_k,f_k)_{k\in \mathbb{N}}\to (\overline{y},\overline{f})$ in $\mathcal{A}$ where $\overline{f}$ is a null control for $\overline{y}$ solution of  \eqref{eq:wave-NL}. Moreover, there exists a $k_0\in \mathbb{N}$ such that the sequence $(\Vert (\overline{y},\overline{f})-(y_k,f_k)\Vert_{\mathcal{A}_0} )_{k\geq k_0}$ decays with a rate equal to $s+1$.
\end{theorem}
\textsc{Proof-} In view of \eqref{estimateF1Y1_A0}, we write 
$$
\sum_{n=0}^k \vert \lambda_n \vert \Vert (Y^1_n,F^1_n)\Vert_{\mathcal{A}_0} \leq C \,m\sum_{n=0}^k  \sqrt{E(y_n,f_n)}.
$$
Using that  
$p_n(\widetilde{\lambda}_n)\leq p_0(\widetilde{\lambda}_0)$ for all $n\geq 0$, we can write for $n>0$,
$$
\sqrt{E(y_n,f_n)}\leq p_{n-1}(\widetilde{\lambda}_{n-1})\sqrt{E(y_{n-1},f_{n-1})}\leq p_{0}(\widetilde{\lambda}_0)\sqrt{E(y_{n-1},f_{n-1})}\leq (p_{0}(\widetilde{\lambda}_0))^n\sqrt{E(y_0,f_0)}.
$$
Then, using that $p_0(\widetilde{\lambda}_{0})=min_{\lambda\in [0,m]}p_0(\lambda)<1$ (since $p_0(0)=1$ and $p_0^\prime(0)<0$), 
we finally obtain the uniform estimate 
$$
\sum_{n=0}^k \vert \lambda_n \vert \Vert (Y^1_n,F^1_n)\Vert_{\mathcal{A}_0} \leq C\, m \frac{\sqrt{E(y_0,f_0)}}{1-p_0(\widetilde{\lambda}_{0})}
$$
for which we deduce that the serie $\sum_{k\geq 0}\lambda_k (Y^1_k,F_k^1)$ converges in $\mathcal{A}_0$. Writing from \eqref{algo_LS_Y} that $(y_{k+1},f_{k+1})=(y_0,f_0)-\sum_{n=0}^k \lambda_n (Y^1_n,F_n^1)$, we conclude that $(y_k,f_k)$ strongly converges in $\mathcal{A}$ to $(\overline{y},\overline{f}):=(y_0,f_0)+\sum_{k\geq 0}\lambda_k (Y^1_k,F_k^1)$. 

Then, using that $(Y^1_k, F^1_k)$ goes to zero as $k\to \infty$ in $\mathcal{A}_0$, we pass to the limit in \eqref{wave-Y1k} and get that $(\overline{y},\overline{f})\in \mathcal{A}$ solves 
\begin{equation}
\label{wave-limit}
\left\{
\begin{aligned}
& \overline{y}_{tt} - \overline{y}_{xx} +  g(\overline{y}) = \overline{f} 1_{\omega}, &\textrm{in}\quad Q_T,\\
& \overline{y}=0, & \textrm{on}\,\, \Sigma_T, \\
& (\overline{y}(\cdot,0),\overline{y}_t(\cdot,0))=(y_0,y_1), & \textrm{in}\,\, \Omega.
\end{aligned}
\right.
\end{equation}
Since the limit $(\overline{y},\overline{f})$ belongs to $\mathcal{A}$, $(\overline{y}(\cdot,T),\overline{y}_t(\cdot,T))=(z_0,z_1)$ in $\Omega$.

We then may write that for all $k>0$
\begin{equation}
\label{estim_coercivity}
\begin{aligned}
\Vert (\overline{y},\overline{f})-(y_k,f_k)\Vert_{\mathcal{A}_0} & =\Vert \sum_{p=k+1}^{\infty} \lambda_p (Y^1_p,F^1_p)\Vert_{\mathcal{A}}\leq m\sum_{p=k+1}^{\infty}  \Vert (Y^1_p,F^1_p \Vert_{\mathcal{A}_0}\\ 
& \leq m\, C\sum_{p=k+1}^{\infty}  \sqrt{E(y_p,f_p)}\\
& \leq m\, C\sum_{p=k+1}^{\infty}  p_{0}(\widetilde{\lambda}_0)^{p-k}\sqrt{E(y_k,f_k)}\\
& \leq m\, C\frac{p_{0}(\widetilde{\lambda}_0)}{1-p_{0}(\widetilde{\lambda}_0)}\sqrt{E(y_{k},f_{k})}
\end{aligned}
\end{equation}
and conclude from Proposition \ref{geomquadratic} the decay with a rate equal to $1+s$ after a finite number of iterates. 
$\hfill\Box$

In particular, along the sequence $(y_k,f_k)_k$ defined by (\ref{algo_LS_Y}), \eqref{estim_coercivity} is a kind of coercivity property for the functional $E$.
We emphasize, in view of the non uniqueness of the zeros of $E$, that an estimate (similar to \eqref{estim_coercivity}) of the form  $\Vert (\overline{y},\overline{f})-(y,f)\Vert_{\mathcal{A}_0} \leq C \sqrt{E(y,f)}$ does not hold for all $(y,f)\in\mathcal{A}$. We also insist in the fact the sequence $(y_k,f_k)_{k>0}$ and its limits $(\overline{y},\overline{f})$ are uniquely determined from the initial guess $(y_0,f_0)$ and from our criterion of selection of the control $F^1$. In other words, the solution $(\overline{y},\overline{f})$ is unique up to the element $(y_0,f_0)\in \A$.


Eventually, if we assume only that $g\in W_0$ and $g^\prime\in L^\infty(\mathbb{R})$, then we get the following result. 
\begin{prop}\label{smallness_gprime}
Assume that $g\in W_0$, $g^\prime\in L^\infty(\mathbb{R})$. If $\Vert g^\prime\Vert_{\infty}$ is small enough so that 
\begin{equation}\label{condition_g_W0}
\sqrt{2}C_1\Vert g^{\prime}\Vert_{\infty}e^{C_2\sqrt{\Vert g^{\prime}\Vert_\infty}}<1.
\end{equation}
then the sequence $(y_k,f_k)_{(k\in \mathbb{N})}$ defined by (\ref{algo_LS_Y}) converges strongly to a controlled pair for \eqref{eq:wave-NL}.
\end{prop}
\textsc{Proof-} Writing, for all $(y,f)\in\mathcal{A}$ and $\lambda\in \mathbb{R}$, that 
$$
\begin{aligned}
2 E((y,f)-\lambda (Y^1,F^1)) & = \Vert (1-\lambda)(y_{tt}-y_{xx}+g(y)-f\,1_{\omega}f)+ g(y-\lambda Y^1)-g(y)+\lambda g^{\prime}(y)Y^1\Vert^2_2\\
& \leq \biggl(\vert 1-\lambda\vert\sqrt{2E(y,f)}+ 2\lambda \Vert g^{\prime}\Vert_{\infty}\Vert Y^1\Vert_2  \biggr)^2.
\end{aligned}
$$
we obtain that 
\begin{equation}
\nonumber
E\big((y_k,f_k)-\lambda (Y^1_k,F_k^1)\big)  \leq   E(y_k,f_k)\biggl(\vert 1-\lambda\vert+\sqrt{2}\lambda\Vert g^{\prime}(y_k)\Vert_{\infty} C_1 e^{C_2\sqrt{\Vert g^{\prime}(y_k)\Vert_{\infty}}}\biggr)^2
\end{equation}
for all $\lambda\in \mathbb{R}$. Taking $\lambda=\lambda_k=1$, the strict decrease of $\sqrt{E(y_k,f_k)}$ w.r.t. $k$ follows if \eqref{condition_g_W0} holds true. 
$\hfill\Box$

\

In the next section, we get similar results of convergence  relaxing the assumption $g^\prime\in L^\infty$.

\subsubsection{Convergence in the case $g\in W_s$, $s\in (0,1]$ and an asymptotic behavior on $g^{\prime}$}

In this section, we assume only that $g^\prime\in L^\infty_{loc}(\mathbb{R})$ and 
\begin{equation}\label{asymptotic_behavior_gprime}
\limsup_{\vert s\vert \to \infty} \frac{\vert g^\prime(s)\vert }{\log^2\vert s\vert} < \beta,
\end{equation}
some constant $\beta>0$. 
Equivalently, we assume that there exists a constant $\alpha>0$ (possibly large) such that 
\begin{equation}\label{asymptotic_behavior_gprime_concrete}
\vert g^{\prime}(s)\vert \leq \alpha + \beta \log^2(1+\vert s\vert), \quad \forall s\in \mathbb{R}.
\end{equation}
The case $\beta=0$ corresponds to the case developed in the previous section, i.e. $g^\prime\in L^\infty(\mathbb{R})$.

Within this more general framework, the difficulty is to have a uniform control with respect to $k$ of the observability constant  $C_1e^{C_2\sqrt{\Vert g^{\prime}(y_k)\Vert_{\infty}}}$ appearing in the estimates for $(Y_k^1,F_k^1)$, see Prop \ref{estimateC1C2}. In other terms, we need to show that the sequence $(y_k,f_k)_{(k\in \mathbb{N})}$ defined in \eqref{algo_LS_Y} is such that $\Vert y_k\Vert_{\infty}$ is uniformly bounded.

\

In the sequel, we define the pair $(y^\star,f^\star)\in \A$ such that $f^\star$ is the control of minimal $L^2(q_T)$-norm for $y^\star$ solution of  \eqref{eq:wave-NL} with $g\equiv 0$.

\begin{lemma}\label{estimate_Ey0f0}
Assume that $g$ satisfies \eqref{asymptotic_behavior_gprime_concrete}. Then, for all $\beta_0>0$, there exists a constant $C(\beta_0)$ such that $\sup_{\beta\in [0,\beta_0]} E(y^\star,f^\star)<C(\beta_0)$.
\end{lemma}
\textsc{Proof}-We get that $E(y^\star,f^\star)=\frac{1}{2}\Vert y_{tt}^\star-y^\star_{xx} +g(y^\star)-f^\star \,1_{\omega}\Vert_2^2=\frac{1}{2}\Vert g(y^\star)\Vert_2^2$. In view of \eqref{asymptotic_behavior_gprime_concrete},
it follows that 
$$
\sqrt{E(y^\star,f^\star)}\leq  T \vert g(0)\vert + T\biggl(\alpha +\beta \log^2(1+\Vert y^\star\Vert_\infty)\biggr)\Vert y^\star\Vert_\infty.
$$
Since $y^\star$ is independent of $g$ and therefore of $\beta$, the result follows.
$\hfill\Box$

\begin{lemma}
Assume that $g$ satisfies \eqref{asymptotic_behavior_gprime_concrete}. For any $(y,f)\in \A$, let $(Y^1,F^1)\in A_0$ be the solution of \eqref{wave-Y1}.
\begin{equation} \nonumber
\Vert (Y^1,Y_t^1)\Vert_{L^{\infty}(0,T;\boldsymbol{V})}+ \Vert F^1\Vert_{2,q_T} \leq d(y)\sqrt{E(y,f)}
\end{equation}
with 
\begin{equation}\nonumber
d(y):=C_3(1+\Vert y\Vert_{\infty})^{C_2 \sqrt{\beta}},  \quad C_3:=C_1 e^{C_2 \sqrt{\alpha}}.
\end{equation}

In particular (up to a constant independent of $y$ and $g$), 
\begin{equation}\label{borneY1k}
\Vert Y\Vert_{L^\infty(Q_T)}\leq d(y)\sqrt{E(y,f)}.
\end{equation}
\end{lemma}
\textsc{Proof-} Using that $\sqrt{a+b}\leq \sqrt{a}+\sqrt{b}$ for all $a,b\geq 0$, it follows that 
$$
e^{C_2\sqrt{\Vert g^{\prime}(y)\Vert_\infty}}\leq  e^{C_2 \sqrt{\alpha}}(1+\Vert y\Vert_{\infty})^{C_2 \sqrt{\beta}}, \quad \forall y\in \A.
$$
\eqref{estimateF1Y1}  then leads to the result.$\hfill\Box$

We also introduce the following notations 
\begin{equation}\label{defc(y)}
c(y):=C \Vert g^{\prime\prime}\Vert_{\infty} \big(d(y)\big)^2, \quad \forall y\in \A
\end{equation}
and 
\begin{defi} Let $\gamma^\star>0$ be defined as follows 
\begin{equation}\nonumber
\gamma^\star=
\left\{
\begin{aligned}
& \frac{-\log\biggl(1-\frac{1}{4c(y^\star)\sqrt{E(y^\star,f^\star)}}\biggr)}{2C_2\log(1+C_3\sqrt{E(y^\star,f^\star)})}, & \textrm{if}\quad 2c(y^\star)\sqrt{E(y^\star,f^\star)}\geq 1,\\
& \frac{1}{2C_2}\frac{\log(2)}{\log(1+C_3\sqrt{E(y^\star,f^\star)})}, & \textrm{else}.
\end{aligned}
\right.
\end{equation}
\end{defi}

\begin{prop}\label{geomquadratic}
Assume $g\in W_s$ for some $s\in (0,1]$ and that $g^\prime$ satisfies the asymptotic behavior \eqref{asymptotic_behavior_gprime} with $\beta<\beta^\star:=\min(C_2^{-2}, (\gamma^\star)^2)$. Let $(y_k,f_k)_{k>0}$ be the sequence of $\mathcal{A}$ defined in \eqref{algo_LS_Y} initialized with $(y_0,f_0)=(y^\star,f^\star)\in \A$. Then $E(y_k,f_k)\to 0$ as $k\to\infty$. Moreover, there exists a $k_0\in \mathbb{N}$ such that the sequence $(E(y_k,f_k))_{k>k_0}$ decays with a rate equal to $s+1$. 
\end{prop}
\textsc{Proof-} In order to simplify the notations, we make the proof in the case $s=1$. \eqref{estimateE} implies that
\begin{equation}\label{estimateEk}
E\big((y,f)-\lambda (Y^1,F^1)\big)\leq   E(y,f)\biggl(\vert 1-\lambda\vert +c(y)\lambda^2 \sqrt{E(y,f)}\biggr)^2.
\end{equation}
 For any $(y_k,f_k)\in \mathcal{A}$, let us then denote the real function $p_k$ by 
$$p_k(\lambda):=\vert 1-\lambda\vert +\lambda^2 c(y_k)\sqrt{E(y_k,f_k)}$$
 for all $\lambda\in [0,1]$. Lemma \ref{lemma_estimateE} with $(y,f)=(y_k,f_k)$ then allows to write that 
$$
\sqrt{E(y_{k+1},f_{k+1})}=\min_{\lambda\in [0,m]} \sqrt{E((y_k,f_k)-\lambda (Y^1_k,F^1_k))}\leq \min_{\lambda\in [0,m]} p_k(\lambda) \sqrt{E(y_k,f_k)}.
$$
We check that the optimal $\lambda_k$ is given by 
\begin{equation}\nonumber
\lambda_k:=\left\{
\begin{aligned}
& \frac{1}{2c(y_k)\sqrt{E(y_k,f_k)}}, & \textrm{if}\quad 2c(y_k)\sqrt{E(y_k,f_k)}\geq 1,\\
& 1, & \textrm{if}\quad 2c(y_k)\sqrt{E(y_k,f_k)}< 1
\end{aligned}
\right.
\end{equation}
leading to 
\begin{equation}
\label{pklambdak}
p_k(\lambda_k):=\left\{
\begin{aligned}
& 1-\frac{1}{4c(y_k)\sqrt{E(y_k,f_k)}}, & \textrm{if}\quad 2c(y_k)\sqrt{E(y_k,f_k)}\geq 1,\\
& c(y_k) \sqrt{E(y_k,f_k)}, & \textrm{if}\quad 2c(y_k)\sqrt{E(y_k,f_k)}< 1.
\end{aligned}
\right.
\end{equation}
Now, assume that $2c(y_k)\sqrt{E(y_k,f_k)}\geq 1$ for some $k\geq 0$. Then, $y_{k+1}=y_k-\lambda_k Y_k^1$ implies that 
$$
\begin{aligned}
(1+\Vert y_{k+1}\Vert_{\infty}) & \leq (1+\Vert y_k\Vert_{\infty})+\lambda_k \Vert Y_k^1\Vert_{\infty},\\
&\leq  (1+\Vert y_k\Vert_{\infty})+ \frac{1}{2c(y_k)\sqrt{E(y_k,f_k)}}d(y_k)\sqrt{E(y_k,f_k)}\\
&\leq  (1+\Vert y_k\Vert_{\infty})+ d(y_k)\sqrt{E(y_k,f_k)}\\
& \leq (1+\Vert y_k\Vert_{\infty})+ C_3(1+\Vert y_k\Vert)^{C_2\sqrt{\beta}}\sqrt{E(y_k,f_k)}\\
&\leq  (1+\Vert y_k\Vert_{\infty})(1+ C_3\sqrt{E(y_k,f_k)})
\end{aligned}
$$
assuming that $C_2\sqrt{\beta}\leq 1$. This leads to 
$$
(1+\Vert y_{k+1}\Vert_{\infty})^{C_2\sqrt{\beta}} \leq (1+\Vert y_k\Vert_{\infty})^{C_2\sqrt{\beta}} \biggl(1+C_3\sqrt{E(y_k,f_k)}\biggr)^{C_2\sqrt{\beta}}.
$$
Consequently, in the case $2c(y_k)\sqrt{E(y_k,f_k)}\geq 1$, we get in view of \eqref{pklambdak}
$$
\sqrt{E(y_{k+1},f_{k+1})}\leq  \biggl(1-\frac{1}{4c(y_k)\sqrt{E(y_k,f_k)}}\biggr) \sqrt{E(y_k,f_k)}.
$$
and therefore 

$$
\begin{aligned}
c(y_{k+1})\sqrt{E(y_{k+1},f_{k+1})}\leq  \biggl(1-\frac{1}{4c(y_k)\sqrt{E(y_k,f_k)}}\biggr) c(y_{k+1})\sqrt{E(y_k,f_k)}.
\end{aligned}
$$
But 
$$
\frac{c(y_{k+1})}{c(y_k)}=\biggl(\frac{d(y_{k+1})}{d(y_k)}\biggr)^2\leq  \biggl(1+C_3\sqrt{E(y_k,f_k)}\biggr)^{2C_2\sqrt{\beta}}.$$
It follows that 
$$
\begin{aligned}
c(y_{k+1})\sqrt{E(y_{k+1},f_{k+1})}\leq  \biggl(1+C_3\sqrt{E(y_k,f_k)}\biggr)^{2C_2\sqrt{\beta}}
\biggl(1-\frac{1}{4c(y_k)\sqrt{E(y_k,f_k)}}\biggr) c(y_k)\sqrt{E(y_k,f_k)}.
\end{aligned}
$$

We start with $k=0$ assuming that the initialization $(y_0,f_0)=(y^\star,f^\star)$ is such that $2 c(y_0)\sqrt{E(y_0,f_0)}\geq~1$. Therefore, 
$$
\begin{aligned}
c(y_1)\sqrt{E(y_1,f_1)}\leq  \biggl(1+C_3\sqrt{E(y_0,f_0)}\biggr)^{2C_2\sqrt{\beta}}
\biggl(1-\frac{1}{4c(y_0)\sqrt{E(y_0,f_0)}}\biggr) c(y_0)\sqrt{E(y_0,f_0)}
\end{aligned}
$$

Now, we take $\beta\in (0,\beta^\star)$ so that 
$$
e(y_0,f_0):=\biggl(1+C_3\sqrt{E(y_0,f_0)}\biggr)^{2C_2\sqrt{\beta}}
\biggl(1-\frac{1}{4c(y_0)\sqrt{E(y_0,f_0)}}\biggr)<1. 
$$
In view of  Lemma \ref{estimate_Ey0f0}, such a $\beta$ strictly positif exists. It follows that 
$c(y_1)\sqrt{E(y_1,f_1)}< c(y_0)\sqrt{E(y_0,f_0)}$.
Repeating the process with $k=1$, still assuming that $2 c(y_1)\sqrt{E(y_1,f_1)}\geq 1$, we have 
$$
c(y_{2})\sqrt{E(y_{2},f_{2})} \leq e(y_1,f_1)\,  c(y_1)\sqrt{E(y_1,f_1)}\leq e(y_0,f_0) c(y_1)\sqrt{E(y_1,f_1)}
$$
using that $e(y_1,f_1)\leq e(y_0,f_0)$ since both $\sqrt{E(y_1,f_1)}<\sqrt{E(y_0,f_0)}$ and $c(y_1)\sqrt{E(y_1,f_1)}<c(y_0)\sqrt{E(y_0,f_0)}$.
It follows that $c(y_2)\sqrt{E(y_2,f_2)}<c(y_1)\sqrt{E(y_1,f_1)}$. Repeating the arguments, we get that the two sequences $(c(y_k)\sqrt{E(y_k,f_k)})_{k>0}$ and $(\sqrt{E(y_k,f_k)})_{k>0}$ strictly decrease. In particular, we get 
$$
c(y_k)\sqrt{E(y_k,f_k)}\leq e(y_0,f_0)^k c(y_0)\sqrt{E(y_0,f_0)}, \quad \forall k\geq 0
$$
as long as $2c(y_k)\sqrt{E(y_k,f_k)}\geq 1$. Its follows that $c(y_k)\sqrt{E(y_k,f_k)}$ and $\sqrt{E(y_k,f_k)}$ goes to zeros. Consequently, there exists a $k_0$ such that $2 c(y_k)\sqrt{E(y_k,f_k)}<1$ for all $k\geq k_0$.

Let $k\geq k_0$ such that $2 c(y_k)\sqrt{E(y_k,f_k)}<1$. Then, the optimal descent step $\lambda_k$ is equal to $1$ and 

$$
\sqrt{E(y_{k+1},f_{k+1})}\leq c(y_k)E(y_k,f_k)=c(y_k)\sqrt{E(y_k,f_k)} \sqrt{E(y_k,f_k)}\leq \frac{1}{2}\sqrt{E(y_k,f_k)}.
$$
Moreover, 
$$
\begin{aligned}
c(y_{k+1})\sqrt{E(y_{k+1},f_{k+1})}& \leq c(y_{k+1})\sqrt{E(y_k,f_k)} \, c(y_k)\sqrt{E(y_k,f_k)}\\
& \leq \biggl(1+C_3\sqrt{E(y_k,f_k)}\biggr)^{2C_2\sqrt{\beta}} \biggl(c(y_k)\sqrt{E(y_k,f_k)}\biggr)\biggl(c(y_k)\sqrt{E(y_k,f_k)}\biggr).
\end{aligned}
$$
But 
$$
\begin{aligned}
\biggl(1+C_3\sqrt{E(y_k,f_k)}\biggr)^{2C_2\sqrt{\beta}} \biggl(c(y_k)\sqrt{E(y_k,f_k)}\biggr)&< \frac{1}{2}\biggl(1+C_3\sqrt{E(y_0,f_0)}\biggr)^{2C_2\sqrt{\beta}}\\
& <\frac{1}{2}\biggl(1-\frac{1}{4c(y_0)\sqrt{E(y_0,f_0)}}\biggr)^{-1}< 1
\end{aligned}
$$
(since  $e(y_0,f_0)<1$) under the assumption that $2c(y_0)\sqrt{E(y_0,f_0)}\geq 1$. This implies that the sequence $(c(y_{k+1})\sqrt{E(y_{k+1})})_k$ decreases strictly and then that the ratio 
$$
\frac{c(y_{k+1})\sqrt{E(y_{k+1},f_{k+1})}}{c(y_{k})\sqrt{E(y_{k},f_{k})}}
$$
decreases as well w.r.t. $k$. It follows that the sequence  $(c(y_k)\sqrt{E(y_k,f_k)})_{k}$ converges to zero as $k$ tends to infinity. 
But, since $c(y_k)=C\Vert g^{\prime\prime}\Vert_{\infty} (d(y_k))^2=C\Vert g^{\prime\prime}\Vert_{\infty}C^2_3(1+\Vert y\Vert_{\infty})^{2C_2 \sqrt{\beta}}\geq C\Vert g^{\prime\prime}\Vert_{\infty}C^2_3$ we get that $E(y_k,f_k)\to 0$ as well.   

Eventually, the relation 
$$
c(y_{k+1})\sqrt{E(y_{k+1},f_{k+1})} 
 \leq \biggl(1+C_3\sqrt{E(y_0,f_0)}\biggr)^{2C_2\sqrt{\beta}} \biggl(c(y_k)\sqrt{E(y_k,f_k)}\biggr)^2
$$
implies the quadratic decrease of $(E(y_k,f_k))_{k\geq k_0}$.

In the more favorable situation for which $2c(y_0)\sqrt{E(y_0)}<1$, we may consider larger values of $\beta$
such that 
$$\biggl(c(y_0)\sqrt{E(y_0,f_0)}\biggr)\biggl(1+C_3\sqrt{E(y_0,f_0)}\biggr)^{2C_2\sqrt{\beta}}<1,$$ 
i.e. 
$$
\sqrt{\beta}<  -\frac{1}{2C_2}\frac{\log(c(y_0)\sqrt{E(y_0,f_0)})}{\log(1+C_3\sqrt{E(y_0,f_0)})}<\frac{1}{2C_2}\frac{\log(2)}{\log(1+C_3\sqrt{E(y_0,f_0)})}.
$$

$\hfill\Box$

 \begin{theorem}\label{Th2.12}
 Assume $g\in W_s$ for some $s\in (0,1]$ and that $g^\prime$ satisfies the asymptotic behavior \eqref{asymptotic_behavior_gprime} with $\beta<\beta^\star:=min(C_2^{-2}, (\gamma^\star)^2)$. Let $(y_k,f_k)_{k>0}$ be the sequence of $\mathcal{A}$ defined in \eqref{algo_LS_Y} initialized with $(y_0,f_0)=(y^\star,f^\star)\in \A$. Then, $(y_k,f_k)_{k\in \mathbb{N}}\to (\overline{y},\overline{f})$ in $\mathcal{A}$ where $\overline{f}$ is a null control for $\overline{y}$ solution of  \eqref{eq:wave-NL}. Moreover, there exists a $k_0\in \mathbb{N}$ such that the sequence $(\Vert (\overline{y},\overline{f})-(y_k,f_k)\Vert_{\mathcal{A}_0} )_{k\geq k_0}$ decays with a rate equal to $s+1$.
\end{theorem}
\textsc{Proof-} The proof is very similar to the proof of Theorem \ref{2.12}. We write that, for any $k>0$,
$$
(y_k,f_k)=(y_0,f_0)+\sum_{n=0}^k \lambda_k (Y_n^1,F_n^1)
$$
leading to, using that $\lambda_k\leq m$ and \eqref{estimateF1Y1_A0},
$
\Vert (y_k,f_k)\Vert_{\A}  \leq \Vert (y_0,f_0)\Vert_{\A} +m\sum_{n=0}^k  d(y_n)\sqrt{E(y_n,f_n)}. 
$
Writing that $d(y)\leq \frac{c(y)}{C \Vert g^{\prime\prime}\Vert_{\infty} C_3 }$, we finally get that 
$$
\begin{aligned}
\Vert (y_k,f_k)\Vert_{\A} &\leq  \Vert (y_0,f_0)\Vert_{\A} +\frac{m}{C\Vert g^{\prime\prime}\Vert_{\infty} C_3 }\sum_{n=0, \atop 2c(y_n)\sqrt{E(y_n,f_n)}\geq 1}^k  c(y_n)\sqrt{E(y_n,f_n)}\\
&\hspace{2cm}+ \frac{m}{C\Vert g^{\prime\prime}\Vert_{\infty} C_3 }\sum_{n=0, \atop 2c(y_n)\sqrt{E(y_n,f_n)}< 1}^k c(y_n)\sqrt{E(y_n,f_n)}.
\end{aligned}
$$

The first sum is finite since the set $I:=\{n\in \mathbb{N},\ 2c(y_n)\sqrt{E(y_n,f_n)}\ge 1\}$ is a finite subset of $\mathbb{N}$. The convergence of the second sum is the consequence of the fact that
$c(y_n)\sqrt{E(y_n,f_n)}$ decays quadratically to zero. 
$\hfill\Box$

\section{Additional comments} \label{sec:remarks}

\textbf{1.} We emphasize that the explicit construction used here allows to recover the null controllability property of \eqref{eq:wave-NL} for nonlinearities $g$ in $W_s$ for one $s\in (0,1]$ satisfying the asymptotic property \eqref{asymptotic_behavior_gprime} on $g^\prime$. Moreover, we do not use a fixed point argument as in \cite{zuazua93}. On the other hand, this asymptotic condition \eqref{asymptotic_behavior_gprime} on $g^\prime$ is slightly stronger than the asymptotic condition \eqref{asymptotic_behavior} made in \cite{zuazua93}: this is due to our linearization of \eqref{eq:wave-NL} which involves $g^{\prime}(s)$ while the linearization \eqref{NL_z} in \cite{zuazua93} involves $g(s)/s$.

Moreover, the additional condition of $g^{\prime}$ in $W_s$, i.e. the existence of one $s\in (0,1)$ such that $\sup_{a,b\in \mathbb{R}, a\neq b} \frac{\vert g^\prime(a)-g^\prime(b)\vert}{\vert a-b\vert^s}< \infty$ allows to get a convergence of the sequence uniformly with respect to the initial guess $(y_0,f_0)\in \A$ and without smallness assumption on the data (see Proposition \ref{smallness_gprime}). In practice, this assumption is not really strong as it suffices to smooth the nonlinear functional $g$. Remark that the functional $g(s)=\alpha+\beta s \log^2(1+\vert s\vert)$, $s\in \mathbb{R}$  for some $\beta>0$ small and any $\alpha$ - which is somehow the limit case in Theorem \ref{zuazuaTH} -  satisfies this assumption (in particular $g^{\prime\prime}\in L^{\infty}(\mathbb{R})$) as well as the asymptotic condition \eqref{asymptotic_behavior_gprime} assumed in this work.

\vspace{0.5cm}
\textbf{2.} 
Among the admissible controlled pair $(y,v)\in \A_0$, we have selected for $(Y_1,F_1)$ solution of \eqref{wave-Y1} the one which minimize the functional $J(v)=\Vert v\Vert^2_{L^2(q_T)}$.This leads to the estimate \eqref{estimateF1Y1} which is the key point in the convergence analysis. The analyze remains true with any other quadratic functional of the form $J(y,v)=\Vert w_1 \, v\Vert^2_{L^2(q_T)} + \Vert w_2 \, y\Vert^2_{L^2(q_T)}$ involving positive weights $w_1$ and $w_2$ (see for instance \cite{cindea_efc_munch_2013}).

\vspace{0.5cm}
\textbf{3.} 
If we introduce $F:\mathcal{A}\to L^2(Q_T)$ by $F(y,f):=(y_{tt}-y_{xx} + g(y)-f\,1_\omega)$, we get that $E(y,f)=\frac{1}{2}\Vert F(y,f)\Vert_{L^2(Q_T)}^2$ and observe that, for $\lambda_k=1$, the algorithm \eqref{algo_LS_Y} coincides with the Newton algorithm associated to the mapping $F$ (mentioned in the introduction, see \ref{Newton-nn}). This explains the super linear convergence of Theorem \ref{Th2.12}, notably a quadratic convergence in the case $s=1$ for which we have a control of $g^{\prime\prime}$ in $L^{\infty}(Q_T)$. The optimization of the parameter $\lambda_k$ allows to get a global convergence of the algorithm and leads to the so-called damped Newton method (for $F$). Under general hypothesis, global convergence for this kind of method is achieved, with a linear rate (for instance; we refer to \cite[Theorem 8.7]{deuflhard}). As far as we know, the analysis of damped type Newton methods  for partial differential equations has deserved very few attention in the literature. We mention \cite{lemoinemunch_time, saramito} in the context of fluids mechanics.   

\vspace{0.5cm}
\textbf{4.} 
Suppose to simplify that $\lambda_k$ equals one (corresponding to the standard Newton method). Then, for each $k$, the optimal pair $(Y_k^1,F_k^1)\in \mathcal{A}_0$ is such that the element $(y_{k+1},f_{k+1})$ minimizes over $\A$ the functional $(z,v)\to J(z-y_k,v-f_k)$ with $J(z,v):=\Vert v\Vert_{L^2(q_T)}$, i.e. the control of minimal $L^2(q_T)$ norm. Instead, we may also select the pair $(Y_k^1,F_k^1)$ such that the element $(y_{k+1},f_{k+1})$ minimizes the functional $(z,v)\to J(z,v)$. This leads to the following sequence $(y_k,f_k)_k$ defined by 
 \begin{equation}
\label{eq:wave_lambdakequal1}
\left\{
\begin{aligned}
& y_{k+1,tt} - y_{k+1,xx} +  g^{\prime}(y_k) y_{k+1} = f_{k+1} 1_{\omega}+g^\prime(y_k)y_k-g(y_k), & \textrm{in}\,\, Q_T,\\
& y_k=0,  & \textrm{on}\,\, \Sigma_T, \\
& (y_{k+1}(\cdot,0), y_{k+1,t}(\cdot,0))=(u_0,u_1), & \textrm{in}\,\, \Omega.
\end{aligned}
\right.
\end{equation}
In this case, for each $k$, $(y_k,f_k)$ is a controlled pair for a linearized wave equation, while, in the case of the algorithm \eqref{algo_LS_Y}, the sequence $(y_k,f_k)$ is a sum of controlled pairs $(Y^1_n,F^1_n)$, $n\leq k$.
This formulation used in \cite{EFC-AM-2012} is different and the corresponding analysis of convergence (at least in the framework of our least-squares setting) is less straightforward because the right hand side term $g^\prime(y_k)y_k-g(y_k)$ is not easily bounded in term  of $\sqrt{E(y_k,f_k)}$.

\vspace{0.5cm}
\textbf{5.} It should be noted as well that the upper bound of the parameter $\beta$ in Theorem \ref{zuazuaTH} depends on $\Omega$ and $T$ but is independent of the initial data $(u_0,u_1)$: precisely,  $\beta< (1+C_2)^{-2}$ where $C_2=C_2(\Omega,T)$ is the constant appearing in \eqref{estimate_u_z}. On the other hand, the upper bound of the parameter $\beta$ in \eqref{asymptotic_behavior_gprime} depends as well on $E(y_0,f_0)$: precisely, 
$\beta<\beta^\star:=\min(C_2^{-2}, (\gamma^\star)^2)$. In particular, for $c(y_0)\sqrt{E(y_0,f_0)}$ small enough, we get that $\beta^\star=C_2^{-2}$ and we recover a bound depending only on $\Omega$ and $T$.

Moreover, as expected, the number  of iterates to achieve convergence (notably to enter in a super-linear regime) depends on the size of the value $c(y_0)\sqrt{E(y_0,f_0)}$. In Theorem \eqref{Th2.12}, we have shown the convergence of the sequence $(y_k,f_k)_{k\in \mathbb{N}}$ when initialized with $(y^\star,f^\star)$ controlled solution of the linear wave equation. This choice is natural and leads to a uniform bound of $E(y_0,f_0)$ in term of $\beta$ in a range $(0,\beta_0)$ (see Lemma \ref{estimate_Ey0f0}).  We may also consider the controlled pair solution of 

\begin{equation}
\label{eq:wave-NL-star-bis}
\left\{
\begin{aligned}
& y^\star_{tt} - y^\star_{xx} +  g(0)+g^\prime(0)y^\star= f^\star 1_{\omega}, & \textrm{in}\,\, Q_T,\\
& y^\star=0, & \textrm{on}\,\, \Sigma_T, \\
& (y^\star(\cdot,0),y^\star_t(\cdot,0))=(u_0,u_1), & \textrm{in}\,\, \Omega,
\end{aligned}
\right.
\end{equation}
leading to $ E(y^\star,f^\star)=\frac{1}{2}\Vert g(y^\star)-g(0)-g^\prime(0)y^\star\Vert^2_{L^2(Q_T)}\leq \frac{1}{2(1+s)^2}\|g^\prime\|^2_{\widetilde{W}^{s,\infty}(\R)}\Vert y^\star\Vert_{\infty}^{2(s+1)}$ and then to $\sqrt{E(y^\star,f^\star)}\leq \|g^\prime\|_{\widetilde{W}^{s,\infty}(\R)}\Vert y^\star\Vert_{\infty}^{(s+1)}$ and also to $\sqrt{E(y^\star,f^\star)}\leq   T\biggl(\alpha +\beta \log^2(1+\Vert y^\star\Vert_\infty)\biggr)\Vert y^\star\Vert_{\infty}$.

\vspace{0.5cm}
\textbf{6.} If the real number $E(y_0,f_0)$ is small enough, then we may remove the asymptotic assumption \eqref{asymptotic_behavior_gprime}  on $g^\prime$.

\begin{prop}
Assume $g\in W_s$ for some $s\in (0,1]$. Let $(y_k,f_k)_{k>0}$ be the sequence of $\mathcal{A}$ defined in \eqref{algo_LS_Y}. There exists a constant $C(\|g^\prime\|_{\widetilde{W}^{s,\infty}(\R)})$ such that if $E(y_0,f_0)\leq C(\|g^\prime\|_{\widetilde{W}^{s,\infty}(\R)})$, then $(y_k,f_k)_{k\in \mathbb{N}}\to (\overline{y},\overline{f})$ in $\mathcal{A}$ where $\overline{f}$ is a null control for $\overline{y}$ solution of  \eqref{eq:wave-NL}. Moreover, there exists a $k_0\in \mathbb{N}$ such that the sequence $(\Vert (\overline{y},\overline{f})-(y_k,f_k)\Vert_{\mathcal{A}_0} )_{k\geq k_0}$ decays with a rate equal to $s+1$.
\end{prop}
\textsc{Proof-} Once again, for simplicity, we make the proof for $s=1$. From \eqref{estimateE}, for all $\lambda\in \mathbb{R}$,
\begin{equation}
\label{ekgk}
\left\{
\begin{aligned}
& \sqrt{E(y_{k+1},f_{k+1})}  \leq   min_{\lambda\in [0,m]}p_k(\lambda)\sqrt{E(y_k,f_k)},\\
& p_k(\lambda):=\vert 1-\lambda\vert + \lambda^2 C(y_k)\sqrt{E(y_k,f_k)},\\
&  C(y_k):=\frac{C}{2\sqrt{2}} \Vert g^{\prime\prime}\Vert_{\infty}  \biggl(C_1  e^{(1+C_2)\sqrt{\Vert g^{\prime}(y)\Vert_\infty}}\biggr)^2.
\end{aligned}
\right.
\end{equation}
We note $B:=\Vert g^{\prime\prime}\Vert_{\infty}$ and $D:=\frac{C}{2\sqrt{2}} B$. 
Then, 
$$
\begin{aligned}
C(y_{k+1}) &=D\,  \biggl(C_1 e^{(1+C_2) \sqrt{\Vert g^{\prime}(y_{k+1})\Vert_\infty}}\biggr)^2= D\,  \biggl(C_1 e^{(1+C_2) \sqrt{\Vert g^{\prime}(y_k-\lambda_k Y_k^1)\Vert_{\infty}}}\biggr)^2\\
&\leq  D  \biggl(C_1 e^{(1+C_2) \sqrt{\Vert g^{\prime}(y_k)\Vert_{\infty}}}e^{(1+C_2)\sqrt{\lambda_k B \Vert Y_k^1\Vert_{\infty}}}\biggr)^2 = C(y_k) \big(e^{(1+C_2)\sqrt{m B}\sqrt{\Vert Y_k^1\Vert_\infty}}\big)^2
\end{aligned}
$$
so that, multiplying \eqref{ekgk} by $C(y_{k+1})$ and introducing the notation $e_{k}:=C(y_k)\sqrt{E(y_k,f_k)}$, we obtain the inequality 
\begin{equation}
e_{k+1}  \leq   \min_{\lambda\in [0,m]}p_k(\lambda)\biggl(e^{(1+C_2)\sqrt{mB} \sqrt{\Vert Y_k^1\Vert_\infty}}\biggr)^2 \, e_k. \label{ekp1ek}
\end{equation}
%
Recalling from \eqref{estimateF1Y1} that $\Vert Y_k^1\Vert_{\infty}\leq C_1\,e^{(1+C_2) \sqrt{\Vert g^{\prime}(y_k)\Vert_{\infty}}}\sqrt{E(y_k,f_k)}=(D^{-1} C(y_k))^{1/2}\sqrt{E(y_k,f_k)}$, we get 
\begin{equation}
(1+C_2)\sqrt{mB} \sqrt{\Vert Y_k^1\Vert_\infty}\leq (1+C_2) \sqrt{m B}  (D^{-1} C(y_k))^{1/4}E(y_k,f_k)^{1/4} \label{ek1/4}
\end{equation}
and therefore 
$$
\biggl(e^{(1+C_2)\sqrt{mB} \sqrt{\Vert Y_k^1\Vert_\infty}}\biggr)^2\leq  e^{C_3 C(y_k)^{1/4} E(y_k,f_k)^{1/4}} 
$$
with $C_3:=2(1+C_2) \sqrt{m B}D^{-1/4}$. Assuming that $E(y_k,f_k)\leq 1$ so that $E(y_k,f_k)^{1/4}\leq E(y_k,f_k)^{1/8}$, we finally get $\big(e^{(1+C_2)\sqrt{mB} \sqrt{\Vert Y_k^1\Vert_\infty}}\big)^2\leq  e^{C_3 e_k^{1/4}}$
and from (\ref{ekp1ek}), 
\begin{equation}
 \nonumber
e_{k+1}  \leq   \min_{\lambda\in [0,m]} \big(\vert 1-\lambda\vert +e_k \lambda^2\big) \quad e^{C_3 e_k^{1/4}}\, e_k.
\end{equation}
If $2e_k<1$, the minimum is reached for $\lambda=1$ leading $\frac{e_{k+1}}{e_k}\leq  e_k e^{C_3 e_k^{1/4}}$. Consequently, if the initial guess $(y_0,f_0)$ belongs to the set $\{(y_0,f_0)\in \mathcal{A}, E(y_0,f_0)\leq 1, e_0< 1/2, e_0 e^{C_3 e_0^{1/4}}<1\}$,
the sequence $(e_k)_{k>0}$ goes to zero as $k\to \infty$ (with a quadratic rate). Since $C(y_k)\geq D C_1^2$ for all $k\in \mathbb{N}$, this implies that the sequence $(E(y_k,f_k))_{k>0}$ goes to zero as well. Moreover, from \eqref{estimateF1Y1_A0}, we get $D \Vert (Y_k^1,F_k^1)\Vert_{\mathcal{A}_0}\leq e_k^{1/2}$ and repeating the arguments of the proof of Theorem \ref{2.12}, we conclude that the sequence $(y_k,f_k)_{k>0}$ converges to a controlled pair for \eqref{eq:wave-NL}.

Remark that these computations are valid under the assumptions that $g^\prime\in L_{loc}^\infty(\mathbb{R})$ and $g^{\prime\prime}\in L^\infty(\mathbb{R})$ but they did not use the assumption \eqref{asymptotic_behavior} nor \eqref{asymptotic_behavior_gprime} on the nonlinearity $g$. 
However, the smallness assumption on $e_0$ requires a smallness assumption on $\sqrt{E(y_0,f_0)}$ (since $d_0>1$). This is equivalent to assume the controllability of \eqref{eq:wave-NL}. Alternatively, in the case $g(0)=0$, the smallness assumption on $\sqrt{E(y_0,f_0}$ is achieved as soon as the data $(u_0,u_1)$  is small for the norm $\boldsymbol{V}$. This result of convergence is therefore equivalent to the local controllability of \eqref{eq:wave-NL}.

\vspace{0.5cm}
\textbf{7.} Under the strong assumption $g^\prime\in L^\infty(\mathbb{R})$ Theorem  \ref{2.12}  remains true in the multi-dimensional case (see \cite{Zhang_2000}) assuming that the triplet $(\Omega,\omega,T)$ satisfies the classical multiplier condition introduced in \cite{JLL88}. 

 \begin{theorem}\label{2.12-ND}
Let $\Omega$ is a bounded subset of $\mathbb{R}^d$, $1\leq d\leq 3$ and $\omega$ a non empty open subset of $\Omega$. Assume that the triplet $(\Omega,\omega,T)$ satisfies the  multiplier condition. 
Assume that $g\in W_s$ for some $s\in (0,1]$ and that $g^\prime\in L^\infty(\mathbb{R})$. Let $(y_k,f_k)_{k\in\mathbb{N}}$ be the sequence defined by (\ref{algo_LS_Y}). Then, $(y_k,f_k)_{k\in \mathbb{N}}\to (\overline{y},\overline{f})$ in $\mathcal{A}$ where $\overline{f}$ is a null control for $\overline{y}$ solution of  \eqref{eq:wave-NL}. Moreover, there exists a $k_0\in \mathbb{N}$ such that the sequence $(\Vert (\overline{y},\overline{f})-(y_k,f_k)\Vert_{\mathcal{A}_0} )_{k\geq k_0}$ decays with a rate equal to $s+1$.
\end{theorem}
The proof, although more technical, follows the line of proof of Theorem \ref{2.12}. We refer to \cite{lemoine_gayte_munch} for the proof in the case of a semi-linear heat equation. Using \cite{Li_Zhang_2000}, Theorem \ref{Th2.12} can be also extended to the multidimensional case replacing the growth condition \eqref{asymptotic_behavior_gprime} by the following one : $\limsup_{\vert s\vert \to \infty} \frac{\vert g^\prime(s)\vert }{\log^{1/2}\vert s\vert} < \beta$ for some $\beta>0$.

\vspace{0.5cm}
\textbf{8.}  Eventually, this approach may be extended, with few modifications, to the boundary case considered notably in \cite{zuazua91}.

\section{Appendix : Controllability results for the linearized wave equation}\label{sec:linearizedwave}

We recall in this appendix some a priori estimates for the linearized wave equation with potential in $L^\infty(Q_T)$ and right hand side in $L^2(Q_T)$. 

\begin{prop}\label{controllability_result}
Let $A\in L^{\infty}(Q_T)$, $B\in L^2(Q_T)$ and $(z_0,z_1)\in \boldsymbol{V}$. Let $\omega=(l_1,l_2)$. Assume $T>2\max(l_1,1-l_2)$. There exist control functions $u\in L^2(q_T)$ such that the solution of 
\begin{equation}
\label{wave_z}
\left\{
\begin{aligned}
& z_{tt} - z_{xx} +  A z  = u 1_{\omega} + B, & \textrm{in}\,\, Q_T,\\
& z=0, & \textrm{on}\,\, \Sigma_T, \\
& (z(\cdot,0),z_t(\cdot,0))=(z_0,z_1), &\textrm{in}\,\, \Omega,
\end{aligned}
\right.
\end{equation}
satisfies $(z(\cdot,T),z_t(\cdot,T))=(0,0)$ in $\Omega$.  Moreover, the unique control $u$ which minimizes the $L^2(q_T)$-norm together with the corresponding controlled solution satisfy the estimate
\begin{equation}\label{estimate_u_z}
\Vert u\Vert_{2,q_T} + \Vert (z,z_t)\Vert_{L^\infty(0,T;\boldsymbol{V})}\leq C_1 \biggl(\Vert B\Vert_2 e^{(1+C_2)\sqrt{\Vert A\Vert_\infty}} + \Vert z_0,z_1\Vert_{\boldsymbol{V}}\biggr) e^{C_2\sqrt{\Vert A\Vert_\infty}}
\end{equation}
for some constant $C_1,C_2>0$. 
\end{prop}

\textsc{Proof -} The proof is based on estimates obtained \cite{zuazua93}. The control of minimal $L^2(q_T)$-norm is given by $u=\ph\, 1_{\omega}$
where $\ph$ solves the adjoint equation 
\begin{equation}
\label{wave_adjoint}
\left\{
\begin{aligned}
& \ph_{tt} - \ph_{xx} +  A \ph  = 0, & \textrm{in}\,\, Q_T,\\
& \ph=0, & \textrm{on}\,\, \Sigma_T, \\
& (\ph(\cdot,0),\ph_t(\cdot,0))=(\ph_0,\ph_1), &\textrm{in}\,\, \Omega,
\end{aligned}
\right.
\end{equation}
with $(\ph_0,\ph_1)\in \boldsymbol{H}:=L^2(\Omega)\times H^{-1}(\Omega)$ the unique minimizer of 
$$
J(\ph_0,\ph_1):=\frac{1}{2}\int\!\!\!\!\int_{q_T} \ph^2 +\int\!\!\!\!\int_{Q_T} B \ph - \langle(z_0,z_1),(\ph_0,\ph_1)\rangle_{\boldsymbol{V},\boldsymbol{H}}
$$
with $\langle(z_0,z_1),(\ph_0,\ph_1)\rangle_{\boldsymbol{V},\boldsymbol{H}}:=\langle z_0,\ph_1\rangle_{H_0^1(\Omega),H^{-1}(\Omega)}-(z_1,\ph_0)_{L^2(\Omega),L^2(\Omega)}$. In particular, the control $v$ satisfies the optimality condition
$$
\int\!\!\!\!\int_{q_T} \ph\,\overline{\ph} +\int\!\!\!\!\int_{Q_T} B \overline{\ph} - \langle(z_0,z_1),(\overline{\ph}_0,\overline{\ph}_1)\rangle_{\boldsymbol{V},\boldsymbol{H}}=0, \quad\forall (\overline{\ph}_0,\overline{\ph}_1)\in \boldsymbol{H}
$$
from which we deduce that $\Vert u\Vert_{2,q_T}^2 \leq \Vert B\Vert_2 \Vert \ph \Vert_2 + \Vert (z_0,z_1)\Vert_{\boldsymbol{V}}\Vert (\ph_0,\ph_1)\Vert_{\boldsymbol{H}}$.
From \cite[Lemma 2]{zuazua93}, we get 
$$
\Vert (\ph,\ph_t)\Vert^2_{L^\infty(0,T;\boldsymbol{H})}\leq B_1 \Vert (\ph_0,\ph_1)\Vert^2_{\boldsymbol{H}} (1+\Vert A\Vert^2_\infty)e^{B_2 \sqrt{\Vert A\Vert_{\infty}}}
$$
for some constant $B_1,B_2>0$ from which it follows that 
$\Vert \ph \Vert^2_2\leq T B_1 \Vert (\ph_0,\ph_1)\Vert^2_{\boldsymbol{H}} (1+\Vert A\Vert^2_{\infty})e^{B_2 \sqrt{\Vert A\Vert_{\infty}}}$.
Moreover, from \cite[Theorem 4]{zuazua93}, there exists $C_1,C_2>0$ such that $\Vert (\ph_0,\ph_1)\Vert^2_{\boldsymbol{H}}\leq C_1e^{C_2 \sqrt{\Vert A\Vert_{\infty}}}\Vert \ph\Vert^2_{2,q_T}$.
Combining these inequalities, we get 
$$
\Vert u\Vert_{L^2(q_T)}\leq \biggl(\Vert B\Vert_2 \sqrt{T}\sqrt{B_1}(1+\Vert A\Vert^2_\infty)^{1/2} e^{\frac{B_2}{2}\sqrt{\Vert A\Vert_\infty}}+\Vert (z_0,z_1)\Vert_{\boldsymbol{V}}\biggr)\sqrt{C_1}e^{\frac{C_2}{2}\sqrt{\Vert A\Vert_{\infty}}}.
$$
Using the inequality $(1+s^2)^{1/2} \leq e^{\sqrt{s}}$ for all $s\geq 0$, we get the result. 
Then, from \cite[Lemma 1]{zuazua93}, we have
$$
\Vert (z,z_t)\Vert^2_{L^\infty(0,T;\boldsymbol{V})}\leq D_1\biggl(  \Vert (z_0,z_1)\Vert^2_{\boldsymbol{H}}(1+\Vert A\Vert_\infty)+ \Vert u\, 1_{\omega}+B\Vert^2_2\biggr)
e^{D_2\sqrt{\Vert A\Vert_\infty}}.
$$
for some constant $D_1,D_2>0$ from which we deduce that 
$$
\Vert (z,z_t)\Vert^2_{L^\infty(0,T;\boldsymbol{V})}\leq D_1\biggl(  \Vert (z_0,z_1)\Vert^2_{\boldsymbol{H}}(2+\Vert A\Vert_\infty)+ 2\Vert B\Vert^2_2\biggl(1+TB_1\big(1+\Vert A\Vert_{\infty}^2\big)e^{B_2\sqrt{\Vert A\Vert_{\infty}}}\biggr) \biggr)
e^{D_2\sqrt{\Vert A\Vert_\infty}}.
$$
Using that $(1+s)^{1/2} \leq e^{\sqrt{s}}$ and  $(1+s^2) \leq e^{2\sqrt{s}}$ for all $s\geq 0$, we get the estimate.  $\hfill\Box$

We then discuss some properties of the operator $K: L^\infty(Q_T)\to L^\infty(Q_T)$ defined by $K(\xi)=y_{\xi}$ a null controlled solution of the linear boundary value problem \eqref{NL_z}
through the control of minimal $L^2(q_T)$ norm $f_{\xi}$. Proposition \ref{controllability_result} with $B=-g(0)$ leads to the estimate
\begin{equation}
\label{estimateK}
\Vert (y_{\xi},y_{\xi,t})\Vert_{L^\infty(0,T;\boldsymbol{V})}\leq C_1\biggl(\Vert u_0,u_1\Vert_{\boldsymbol{V}}+\Vert g(0)\Vert_2 e^{(1+C_2)\sqrt{\Vert \widehat{g}(\xi)\Vert_{\infty}}} \biggr)e^{C_2\sqrt{\Vert \widehat{g}(\xi)\Vert_\infty}}.
\end{equation}
Then, as in \cite{zuazua93}, we write that the assumption \eqref{assumption_g} on $g$ implies that there exists some $d>0$ such that $\Vert \widehat{g}(y)\Vert_\infty\leq d+\beta \log^2(1+\Vert y\Vert_\infty)$
for all $y\in L^\infty(Q_T)$ from which it follows that $e^{C_2\sqrt{\Vert \hat{g}(\xi)\Vert_\infty}}\leq e^{C_2\sqrt{d}} (1+\Vert \xi\Vert_{\infty})^{C_2\sqrt{\beta}}$.
\eqref{estimateK} then leads to the estimate
$$
 \Vert y_{\xi}\Vert_{\infty}\leq C_1\biggl(\Vert u_0,u_1\Vert_{\boldsymbol{V}}+\Vert g(0)\Vert_2\biggr)e^{(1+2C_2)\sqrt{d}} (1+\Vert \xi\Vert_{\infty})^{(1+2C_2)\sqrt{\beta}}.
$$
Taking $\beta$ small enough so that $(1+2C_2)\sqrt{\beta}<1$, we conclude that there exists a constant $M>0$ such that $\Vert \xi\Vert_\infty\leq M$ implies $\Vert K(\xi)\Vert_\infty\leq M$. This is the argument in \cite{zuazua93}. Contrary to $\beta$, we remark that $M$ depends on $\Vert u_0,u_1\Vert_{\boldsymbol{V}}$ (and increases with $\Vert u_0,u_1\Vert_{\boldsymbol{V}}$).

The following proposition gives an estimate of the difference of two controlled solutions. 

\begin{prop}\label{gap}
Let $a,A\in L^\infty(Q_T)$ and $B\in L^2(Q_T)$. Let $u$ and $v$ be the null controls of minimal $L^2(q_T)$-norm for $y$ and $z$ solutions of 
\begin{equation}
\label{wave_y}
\left\{
\begin{aligned}
& y_{tt} - y_{xx} +  A y  = u 1_{\omega} + B, & \textrm{in}\,\, Q_T,\\
& y=0, & \textrm{on}\,\, \Sigma_T, \\
& (y(\cdot,0),y_t(\cdot,0))=(u_0,u_1), &\textrm{in}\,\, \Omega,
\end{aligned}
\right.
\end{equation}
and 
\begin{equation}
\label{wave_zz}
\left\{
\begin{aligned}
& z_{tt} - z_{xx} +  (A+a)z   = v 1_{\omega} + B, & \textrm{in}\,\, Q_T,\\
& z=0, & \textrm{on}\,\, \Sigma_T, \\
& (z(\cdot,0),z_t(\cdot,0))=(u_0,u_1), &\textrm{in}\,\, \Omega,
\end{aligned}
\right.
\end{equation}
respectively. Then, 
$$
\begin{aligned}
\Vert y-z\Vert_{L^\infty(Q_T)} 
& \leq C_1^4 \Vert a \Vert_\infty e^{C_2\sqrt{\Vert A+a\Vert_\infty}} e^{(2+3C_2)\sqrt{\Vert A\Vert_\infty}}\biggl(\Vert B\Vert_2 e^{(1+C_2)\sqrt{\Vert A\Vert_\infty}} + \Vert u_0,u_1\Vert_{\boldsymbol{V}}\biggr).
\end{aligned}
$$
for some constants $C_1,C_2>0$.
\end{prop}
\textsc{Proof}- We write that the control of minimal $L^2$-norm for $y$ and $z$ are given by $u=\ph\, 1_\omega$ and 
$v=\ph_a\, 1_\omega$ where $\ph$ and $\ph_a$ solve the adjoint equation  
 \begin{equation}
 \nonumber
\left\{
\begin{aligned}
& \ph_{tt} - \ph_{xx} +  A \ph   = 0, & \textrm{in}\,\, Q_T,\\
& \ph=0, & \textrm{on}\,\, \Sigma_T, \\
& (\ph(\cdot,0),\ph_t(\cdot,0))=(\ph_0,\ph_1), &\textrm{in}\,\, \Omega,
\end{aligned}
\right.
\qquad\quad 
\left\{
\begin{aligned}
& \ph_{a,tt} - \ph_{a,xx} +  (A+a) \ph_a   = 0, & \textrm{in}\,\, Q_T,\\
& \ph=0, & \textrm{on}\,\, \Sigma_T, \\
& (\ph(\cdot,0),\ph_t(\cdot,0))=(\ph_{a,0},\ph_{a,1}), &\textrm{in}\,\, \Omega,
\end{aligned}
\right.
\end{equation}
for some appropriate $(\ph_0,\ph_1), (\ph_{a,0},\ph_{a,1})\in \boldsymbol{H}$. Consequently, the difference $Z:=z-y$ solves 
\begin{equation}
\label{diffZ}
\left\{
\begin{aligned}
& Z_{tt} - Z_{xx} +  (A+a)Z   = \Phi 1_{\omega} -ay, & \textrm{in}\,\, Q_T,\\
& Z=0, & \textrm{on}\,\, \Sigma_T, \\
& (Z(\cdot,0), z_t(\cdot,0))=(0,0), &\textrm{in}\,\, \Omega,
\end{aligned}
\right.
\end{equation}
while $\Phi:=(\ph_a-\ph)$ solves
\begin{equation}
 \nonumber
\left\{
\begin{aligned}
& \Phi_{tt} - \Phi_{xx} +  (A+a)\Phi  = -a\ph & \textrm{in}\,\, Q_T,\\
& \Phi=0, & \textrm{on}\,\, \Sigma_T, \\
& (\Phi(\cdot,0), \Phi_t(\cdot,0))=(\ph_{a,0}-\ph_0,\ph_{a,1}-\ph_1), &\textrm{in}\,\, \Omega.
\end{aligned}
\right.
\end{equation}
We decompose $\Phi=\Psi+\psi$ where $\Psi$ and $\psi$ solves respectively 
\begin{equation}
 \nonumber
\left\{
\begin{aligned}
& \Psi_{tt} - \Psi_{xx} +  (A+a)\Psi  = 0 & \textrm{in}\,\, Q_T,\\
& \Psi=0, & \textrm{on}\,\, \Sigma_T, \\
& (\Psi(\cdot,0), \Psi_t(\cdot,0))=(\ph_{a,0}-\ph_0,\ph_{a,1}-\ph_1), &\textrm{in}\,\, \Omega,
\end{aligned}
\right.
\qquad \left\{
\begin{aligned}
& \psi_{tt} - \psi_{xx} +  (A+a)\psi   =  -a\ph & \textrm{in}\,\, Q_T,\\
& \psi=0, & \textrm{on}\,\, \Sigma_T, \\
& \psi(\cdot,0), \psi_t(\cdot,0))=(0,0), &\textrm{in}\,\, \Omega,
\end{aligned}
\right.
\end{equation}
and then deduce that $\Psi\,1_\omega$ is the control of minimal $L^2$-norm for $Z$ solution of  
\begin{equation}
 \nonumber
\left\{
\begin{aligned}
& Z_{tt} - Z_{xx} +  (A+a)Z   = \Psi 1_{\omega} + \biggl(\psi 1_\omega-ay\biggr), & \textrm{in}\,\, Q_T,\\
& Z=0, & \textrm{on}\,\, \Sigma_T, \\
& (Z(\cdot,0), Z_t(\cdot,0))=(0,0), &\textrm{in}\,\, \Omega.
\end{aligned}
\right.
\end{equation}
Proposition \ref{controllability_result} implies that 

\begin{equation}\nonumber
\Vert \Psi\Vert_{2,q_T}+ \Vert (Z,Z_t)\Vert_{L^\infty(0,T;\boldsymbol{V})}\leq C_1 \Vert \psi 1_\omega-ay\Vert_2 e^{(1+2C_2)\sqrt{\Vert A\Vert_\infty}}.
\end{equation}
Moreover, energy estimates for $\psi$ leads to 
$
\Vert \psi\Vert_{L^2(q_T)}\leq C_1\Vert a\Vert_\infty \Vert\ph\Vert _{2}e^{C_2\sqrt{\Vert A+a\Vert_\infty}}
$
and also to 
$$
\Vert\ph\Vert _2\leq C_1 \Vert \ph_0,\ph_1\Vert_{\boldsymbol{H}}e^{(1+C_2)\sqrt{\Vert A\Vert_\infty}}\leq \biggl(C_1 e^{(1+C_2)\sqrt{\Vert A\Vert_\infty}}\biggr)^2 \Vert u\Vert_{2,q_T}
$$
using that $\Vert \ph_0,\ph_1\Vert_{\boldsymbol{H}}\leq C_1 e^{C_2\sqrt{\Vert A\Vert_\infty}}\Vert u\Vert_{2,q_T}$
so that 
$$
\Vert \psi\Vert_{L^2(q_T)}\leq C_1\Vert a\Vert_\infty e^{C_2\sqrt{\Vert A+a\Vert_\infty}} \biggl(C_1 e^{(1+C_2)\sqrt{\Vert A\Vert_\infty}}\biggr)^2 \Vert u\Vert_{2,q_T}
$$
from which we deduce that
$$
\begin{aligned}
\Vert Z\Vert_{L^\infty(Q_T)} & \leq C_1 \biggl(\Vert \psi 1_\omega\Vert + \Vert a\Vert_{L^\infty(Q_T)}\Vert y\Vert_2\biggr) e^{(1+2C_2)\sqrt{\Vert A\Vert_\infty}}\\
& \leq C_1 \Vert a \Vert_\infty \biggl(e^{C_2\sqrt{\Vert A+a\Vert_\infty}} \biggl(C_1 e^{(1+C_2)\sqrt{\Vert A\Vert_\infty}}\biggr)^2 \Vert u\Vert_{L^2(q_T)}+\Vert y\Vert_{L^2(Q_T)}\biggr)\\
& \leq C_1 \Vert a \Vert_\infty \biggl(e^{C_2\sqrt{\Vert A+a\Vert_\infty}} \biggl(C_1 e^{(1+C_2)\sqrt{\Vert A\Vert_\infty}}\biggr)^2+1\biggr)\\
&\hspace{2cm}C_1 \biggl(\Vert B\Vert_2 e^{(1+C_2)\sqrt{\Vert A\Vert_\infty}} + \Vert u_0,u_1\Vert_{\boldsymbol{V}}\biggr) e^{C_2\sqrt{\Vert A\Vert_\infty}}\\
\end{aligned}
$$
leading to the result. 
$\hfill\Box$

This result allows to show the following property on the operator $K$.

\begin{prop}
Assume hypotheses of Theorem \ref{zuazuaTH}.
Let $M=M(\Vert u_0, u_1\Vert_{\boldsymbol{V}},\beta)$ a constant such that $K$ maps $B_\infty(0,M)$ into itself and assume that $\hat{g}^{\prime}\in L^\infty(0,M)$.
For any $\xi^i\in B_\infty(0,M)$, $i=1,2$, there exists a constant $c(M)>0$ such that  
\begin{equation}
 \nonumber
\Vert K(\xi^2)-K(\xi^1)\Vert_{\infty}\leq c(M) \Vert \hat{g}^{\prime}\Vert_{L^\infty(0,M)} \Vert \xi^2-\xi^1\Vert_\infty.
\end{equation}
\end{prop}
\textsc{Proof-} For any $\xi^i\in B_\infty(0,M)$, $i=1,2$, let $y_{\xi^i}=K(\xi^i)$ be the null controlled solution of 
\begin{equation}
 \nonumber
\left\{
\begin{aligned}
& y_{\xi^i,tt} - y_{\xi^i,xx} +  y_{\xi^i} \,\widehat{g}(\xi^i)= -g(0)+f_{\xi^i} 1_{\omega}, &\textrm{in}\,\, Q_T,\\
& y_\xi^i=0, &\textrm{on}\,\, \Sigma_T, \\
& (y_{\xi^i}(\cdot,0),y_{\xi^i,t}(\cdot,0))=(u_0,u_1), &\textrm{in}\,\, \Omega,
\end{aligned}
\right.
\qquad
\widehat{g}(s):=
\left\{ 
\begin{aligned}
& \frac{g(s)-g(0)}{s} & s\neq 0,\\
& g^{\prime}(0) & s=0
\end{aligned}
\right., 
\end{equation}
through the control of minimal $L^2(q_T)$ norm $f_{\xi^i} 1_{\omega}$. We observe that $y_{\xi^2}$ is solution of 
\begin{equation}
 \nonumber
\left\{
\begin{aligned}
& y_{\xi^2,tt} - y_{\xi^2,xx} +  y_{\xi^2} \,\widehat{g}(\xi^1)+ y_{\xi^2}(\widehat{g}(\xi^2)-\widehat{g}(\xi^1))= -g(0)+f_{\xi^2} 1_{\omega}, &\textrm{in}\,\, Q_T,\\
& y_{\xi^2}=0, &\textrm{on}\,\, \Sigma_T, \\
& (y_{\xi^2}(\cdot,0),y_{\xi^2,t}(\cdot,0))=(u_0,u_1), &\textrm{in}\,\, \Omega.
\end{aligned}
\right.
\end{equation}
Therefore, from Proposition \ref{gap} with $B=-g(0)$, $A=\hat{g}(\xi^1)$, $a=\hat{g}(\xi^2)-\hat{g}(\xi^1)$, $y_{\xi^2}-y_{\xi^1}$ satisfies 
\begin{equation}
\label{gapxi}
\Vert y_{\xi^2}-y_{\xi^1} \Vert_{\infty}\leq A(\xi^1,\xi^2)\Vert \widehat{g}(\xi^2)-\widehat{g}(\xi^1)\Vert_\infty
\end{equation}
where the positive constant  
$$
\begin{aligned}
A(\xi^1,\xi^2):=& C_1^2 \biggl(e^{C_2\sqrt{\Vert  \hat{g}(\xi^2)\Vert_\infty}} \biggl(C_1 e^{(1+C_2)\sqrt{\Vert  \hat{g}(\xi^1)\Vert_\infty}}\biggr)^2\biggr)\\
&\biggl(\Vert g(0)\Vert_2 e^{(1+C_2)\sqrt{\Vert  \hat{g}(\xi^1)\Vert_\infty}} + \Vert u_0,u_1\Vert_{\boldsymbol{V}}\biggr) e^{C_2\sqrt{\Vert  \hat{g}(\xi^1)\Vert_\infty}}
\end{aligned}
$$

is bounded in term of $M$ for all $\xi^i\in B_\infty(0,M)$. We introduce $c(M)$ such that $A(\xi^1,\xi^2)\leq c(M)$.  \eqref{gapxi} then leads to the result. 

In particular, if $\Vert \hat{g}^{\prime}\Vert_{L^\infty(0,M)}\leq \frac{d}{C(M)}$ for some $d<1$, then the operator $K$ is contracting. Remark however that the bound depends on the norm $\Vert u_0,u_1\Vert_{\boldsymbol{V}}$ of the initial data to be controlled.
$\hfill\Box$

\par\noindent

\vskip 0.5cm
\par\noindent
\textbf{Acknowlegment-}
The first author warmly thanks J\'erome Lemoine (Laboratoire Math\'ematique Blaise Pascal, Clermont Auvergne University) for fruitful discussions concerning this work.

 \bibliographystyle{siam}
 \bibliography{wavecontinuation.bib}
\end{document}